\documentclass{article}

\usepackage{amsfonts}
\usepackage{amssymb}
\usepackage{amsmath}
\usepackage{amsthm}
\usepackage{latexsym}

\newtheorem{thm}{Theorem}[section]
\newtheorem{lemma}[thm]{Lemma}
\newtheorem{prop}[thm]{Proposition}
\newtheorem{cor}[thm]{Corollary}

\newtheorem{rem}[thm]{Remark}

\theoremstyle{definition}

\newcommand{\lp}{\left(}
\newcommand{\rp}{\right)}
\newcommand{\texto}[1]{\quad\mbox{#1}\quad}

\newcommand{\abs}[1]{\left\lvert #1\right\rvert}
\newcommand{\norm}[1]{\left\lVert #1\right\rVert}
\newcommand{\pd}[2]{\frac{\partial #1}{\partial #2}}
\newcommand{\spd}[2]{\frac{\partial^2 #1}{\partial {#2}^2}}

\newcommand{\mR}{\mathbb{R}}

\newcommand{\V}{\mathcal V}
\newcommand{\T}{\mathbb T^1}

\newcommand{\be}{\begin{equation}}
\newcommand{\ee}{\end{equation}}
\newcommand{\bee}{\begin{equation*}}
\newcommand{\eee}{\end{equation*}}
\newcommand{\bea}{\begin{eqnarray}}
\newcommand{\eea}{\end{eqnarray}}
\newcommand{\bs}{\begin{split}}
\newcommand{\es}{\end{split}}

\begin{document}

\title{{\bf Asymptotics for a free-boundary model \\in price formation}}

 \author{Mar\'ia del Mar Gonz\'alez \thanks{Supported by Spain Government project MTM2008-06349-C03-01,
GenCat 2009SGR345 and NSF grant DMS-0635607.}\\Univ. Polit\`ecnica de Catalunya  \and Maria Pia Gualdani\\UT Austin \thanks{Supported by the NSF Grant DMS-0807636.}}

\date{}

\maketitle

\begin{abstract}
We study the asymptotics for large time of solutions to a one dimensional parabolic evolution equation with
non-standard measure-valued right hand side, that involves derivatives of
the solution computed at a free boundary point. The problem is a particular case of a mean-field free boundary model proposed by Lasry-Lions on price formation and dynamic equilibria.

The main step in the proof is based on the fact that the free boundary disappears in the linearized problem, thus can be treated as a perturbation through semigroup theory. This requires a delicate choice for the function spaces since higher regularity is needed near the free boundary. We show global existence for solutions with initial data in a small neighborhood of any equilibrium point, and exponential decay towards a stationary state. Moreover, the family of equilibria of the equation is stable, as follows from center manifold theory.
\end{abstract}


2000 {\bf MSC:} 35R35, 35K15, 91B42, 91B26

%------------------------------------------------------------------------------%-------------------------
%
%-------------------------------------------------------------------------------------------------------

\section{Introduction}

\setcounter{equation}{00}

We consider an idealized population of players consisting of
two groups, namely one group of buyers of a certain good and one group of
vendors of the same good. The two groups are described by two non-negative densities
$f_B$, $f_V$ depending on $(x,t)\in\mR\times\mR_+$. In the model, $x$ denotes a
possible value of the price and $t$ the time.

At a certain time $t$, the vendors would like to sell the good, and the function $f_V(x,t)$ describes the density of the vendors that are wiling to sell the good at price $x$. Meanwhile the buyers will
try to get the good at a cheaper price.  The transaction takes place when the two groups agree on the price: we denote by $p(t)$ the {\em{agreement}} price. The price $p(t)$ will be the highest price the buyers are willing to pay, and the lowest price the vendors agreed to sell the good. There exists a transaction cost, which is denoted by a positive constant $a$. When a buyer gets the good for the price $p(t)$, the actual cost of his trade is $p(t) +a$, as well as the profit for the seller is $p(t) -a$. As a consequence, the buyer that got the good for the price $p(t)$, will try in a later time to sell the good at least at the price $p(t) +a$ and the vendor that sold the good for $p(t)$ will try to get at a later time the same good for a price not higher than $p(t) -a$. Thus the parameter $a$ introduces some friction in the system.

The price $p(t)$ results from a dynamical equilibrium between the two density functions. The randomness in the problem is measured by the diffusion coefficient of the two densities $f_V$ and $f_B$, and is denoted by a parameter $\sigma>0$.

The above situation can be described by
the following system of free-boundary evolution equations:
\begin{equation}
\label{equation-buyers}\left\{ \begin{array}{ll}
&\pd {f_B} t -\frac{\sigma^2}{2}\spd{f_B}{x}=\lambda(t)\delta_{x=p(t)-a}
\quad \hbox{if}\: x\le p(t), \; t>0,\\
&f_B(x,t) > 0 \textrm{ if } x<p(t), \quad f_B(x,t)=0 \textrm{ if } x\geq p(t),
\end{array}\right.
\end{equation}
together with
\be
\label{equation-vendors}\left\{ \begin{array}{ll}
&\pd {f_V} t -\frac{\sigma^2}{2}\spd{f_V}{x}=\lambda(t)\delta_{x=p(t)+a} \quad
\hbox{if } x>p(t), \; t>0,\\
&f_V(x,t) >0 \quad \textrm{if } x>p(t), \quad f_V(x,t)=0 \textrm{ if } x\leq p(t),
\end{array}\right.
\ee
where
\be\label{equality-flux}
\lambda(t)=-\frac{\sigma^2}{2}\pd {f_B}x \lp p(t),t\rp=\frac{\sigma^2}{2}\pd
{f_V}x \lp p(t),t\rp.\ee
The symbol $\delta$ denotes the Dirac delta at the indicated
point. The multiplier $\lambda(t)$ represents the number of transactions
at time $t$, so \eqref{equality-flux} means that the flux of buyers which must
be equal to the flux of
vendors. The initial conditions
\bee
f_B(x,0)=f_B^I(x) \mbox{ and }f_V(x,0)=f_V^I(x)
\eee
are such that, for some $p_I$ in $\mR$,
\bee
\begin{split}
&f_B^I(x)>0 \hbox{ if } x<p_I, \quad f_B^I(x)=0 \hbox{ if } x\geq p_I \\
&f_V^I(x)>0 \hbox{ if } x>p_I, \quad f_B^I(x)=0 \hbox{ if } x\leq p_I.
\end{split}
\eee
The equation satisfy the property of conservation of mass. Indeed, both
$$\int_{-\infty}^{p(t)} f_V \,dx\texto{and}\int_{p(t)}^{+\infty} f_B\,dx$$
remain constant for all $t\geq 0$.\\

Equations \eqref{equation-buyers}-\eqref{equality-flux} describe a mean-field model for the
dynamical formation of the price of some good that has been very recently introduced in
\cite{Lasry-Lions}. There the authors  proved that, under suitable assumptions on
smoothness and integrability on the initial data,
there exists a unique smooth solution $(f(x,t),p(t))$ for  $x\in \mathbb R$ of
\eqref{equation-buyers}-\eqref{equality-flux}.

The most important question we are going to address here concerns the long time behavior of the system: will the good reach a stable price ( $p(t) \to \textrm{constant}$ as $t\to \infty$ ?) or will the price keep oscillating in time and never reach a stable value?\\

We remark here that in a bounded interval with symmetric initial data then the solution remains symmetric for all times and asymptotics were proved by the authors in their previous work \cite{price-formation}. However, the general case contains a new ingredient: a free boundary (see \cite{Caffarelli-Salsa}, for instance, for some background and examples on these type of problems).\\

In this work we address the problem  \eqref{equation-buyers}-\eqref{equality-flux} in a bounded interval $[-A,B]$, $A,B>0$, for $a<\min\{A/2,B/2\}$, with zero-Neumann boundary conditions. The aim is to show that if we start with an initial condition that is near a general equilibrium point in some suitable function space, then there exists a unique solution of  \eqref{equation-buyers}-\eqref{equality-flux}, that decays exponentially fast in time to a unique stationary state. In addition, in can be shown that the problem presents a two-dimensional family of equilibria, and that this family is stable.

Although there is a well developed theory of semigroups and invariant manifolds for the study of evolution equations, the main novelty here is the fact that dynamical system arguments can be used for a problem that presents a free boundary. This is possible since we succeeded to treat the free boundary as a perturbation of the linearized problem. In fact, as we will see in the following sections, the free boundary disappears in the linearization and appears again in the nonlinear part of the problem as a term of lower order.

This allows us to study the linearized operator with the classical semigroup theory and to get time estimates for the linear equation. Even though the linearized operator is a non standard one, we can explicitly compute its eigenvalues and corresponding eigenfunctions. Unfortunately the eigenfunctions do not build an orthogonal basis with respect to the standard product in $L^2$. This complicates the choice of functional spaces.

Indeed, the choice of function spaces is a delicate step in the proof. They need to be big enough to allow delta functions in the equation, but on the other hand, higher regularity is needed near the free boundary. In order to give a more explicit characterization of those spaces, interpolation theory and pseudo-differential operators are needed. \\

For simplicity of the notation, we rewrite the problem \eqref{equation-buyers}-\eqref{equality-flux} as the single equation
\be\label{equation-introduction}
\left\{\begin{split}
\pd{f}{ t} -\frac{\sigma^2}{2}\frac{\partial^2 f}{\partial x^2}
&=\lambda(t)\left[ \delta_{x=p(t)-a}-\delta_{x=p(t)+a}\right]\quad \mbox{in }[-A,B]\times \mathbb R_+, \\
f(x,0)&=f_I(x) \quad \mbox{in }[-A,B],\\
f_x(-A,0)& =f_x(B,0)=0,
\end{split}\right.
\ee
where
\bee\lambda(t):=-\frac{\sigma^2}{2} \pd f x(p(t),t),\eee
and
\bee
f:= f_B - f_V, \; f_I = f_B^I - f_V^I, \;\textrm{and}\; p(0)=p_I.
\eee
Moreover, for $t>0$,
$$f(p(t),t)=0,\quad f(x,t)> 0\mbox{ for all } x\in [-A,p(t)),\quad f(x,t)< 0\mbox{ for all } x\in (p(t),B].$$
The unknowns are two quantities: $f(x,t)$ (the solution) and $p(t)$ (the free boundary).
Assume, without loss of generality, that $\sigma^2/2=1$.\\

The problem satisfies important conservation laws. Indeed, condition \eqref{equality-flux} implies preservation of goods and players:
$$\int_{-A}^{p(t)} f \,dx=m_1\texto{and}-\int_{p(t)}^{B} f\,dx=m_2,\quad \mbox{for all } t>0,$$
with $m_1:=\int_{-A}^{p_I} f_I \,dx$ and $m_2:=-\int_{p_I}^{B} f_I \,dx$, $m_1,m_2>0$. Note that $m_1+ m_2$ represents the total amount of players and $m_2$ the amount of goods.

\bigskip

There is an equivalent way to write the problem as a coupled system of equations.  To see it,
we differentiate  the equation $f(p(t),t)=0$, to obtain that
\be
p'(t)=-\frac{f_{xx}(p(t),t))}{f_x(p(t),t)}.
\ee
Thus our problem is equivalent to the system
\be\label{systemPDE}
\left\{\begin{split}
f_t& =f_{xx}-f_x(p(t),t) \left[\delta_{x=p(t)-a}-\delta_{x=p(t)+a}\right],\\
p'(t) &=-\frac{f_{xx}(p(t),t))}{f_x(p(t),t)},
\end{split}\right.\ee
with boundary conditions
$$f_x(-A,t)=f_x(B,t)=0,$$
and initial data
$$f(x,0)=f_I(x), \;p(0)=p_I.$$

\bigskip

To finish, let us remark that the present work (existence and decay results for small initial data) is the first step in the study of general asymptotics for this problem, which is an ongoing project of the authors with L. Chayes and I. Kim. In particular, in the recent paper \cite{Chayes-Gonzalez-Gualdani-Kim}, global existence in time and uniqueness for general initial data has been shown. The present paper is the only available information so far on the asymptotic behavior.

Note that Chayes-Kim have recently studied an evolution Stefan problem in \cite{Chayes-Kim:two-sided}, \cite{Chayes-Kim:supercooled} in relation to particle systems. Although apparently unrelated to our problem, it shares many of its features.

% ==============================================================================================
% ==============================================================================================

\section{Main results}

\setcounter{equation}{00}

First we set up the problem in the interval $[-A,B]$, for $A,B>0$: we seek functions $f(x,t)$ and $p(t)\in(-A+a,B-a)$, $t\in[0,+\infty)$ that solve
\be\label{non-linear}
\left\{\begin{split}
&f_t =f_{xx}-f_x(p(t),t) \left[\delta_{x=p(t)-a}-\delta_{x=p(t)+a}\right],\\
&f(p(t),t)=0,\quad f(x,t)> 0\mbox{ for }x\in[-A,p(t)),\quad f(x,t)< 0\mbox{ for }x\in(p(t),B],\\
&f_x(-A, t)=f_x(B,t)=0,
\end{split}\right.\ee
with initial data $f_I(x)$, $p_I$ satisfying
\be\label{initial-data}
\left\{\begin{split}
&f_I(p_I)=0,\\
&f_I(x)> 0\mbox{ for }x\in[-A,p_I),\quad f_I(x)< 0\mbox{ for }x\in(p_I,B],\\
&f'_I(-A)=f'_I(B)=0,\\
& p(0)=p_I.
\end{split}\right.\ee
The initial datum has mass
\begin{align}\label{mass-init-data}
\int_{-A}^{p_I} f_I \,dx=m_1, \quad \quad -\int_{p_I}^{B} f_I\,dx=m_2,\quad m_1,m_2>0.
\end{align}

The equilibrium states of the above problem are well understood in Section \ref{section-stationary}. First, note that system \eqref{non-linear}-\eqref{initial-data} has infinitely many equilibrium points, i.e., functions which satisfy equation
\begin{align}\label{stationary-equation}
f^0_{xx}-f^0_x(p^0) \left[\delta_{x=p^0-a}-\delta_{x=p^0+a}\right]=0,
\end{align}
where $f^0(p^0)=0$. More precisely
\be\label{equilibrium} f^0(x)=\left\{\begin{split}
&-\lambda^0(x-p^0), & x\in(p^0-a,p^0+a), \\
&-\lambda^0 a, &x \in(p^0+a,B),\\
&\lambda^0a, &\quad  x\in (-A,p^0-a),
\end{split}\right.
\ee
for $\lambda^0 >0$ and $p^0 \in ( -A+a, B-a)$.\\

We denote with $f_\infty$ the {\em{unique}} equilibrium point that is solution of the stationary problem  (\ref{non-linear})-(\ref{mass-init-data}) (which means the only equilibrium point $f^0$ that satisfies also the preservation of mass condition), and it is given by
\bee f_\infty(x)=\left\{\begin{split}
&-\lambda_\infty(x-p_\infty), & x\in(p_\infty-a,p_\infty+a), \\
&-\lambda_\infty a, &x \in(p_\infty+a,B),\\
&+\lambda_\infty a, &\quad  x\in (-A,p_\infty-a),
\end{split}\right.
\eee
where $p_\infty$ and $\lambda_\infty$ are uniquely determined by the conservation of mass condition
$$
 \int_{-A}^{p_\infty} f_\infty \,dx=m_1, \quad \quad -\int_{p_\infty}^{B} f_\infty\,dx=m_2,
 $$
and $m_1$ and $m_2$ are the constants given by the initial datum in \eqref{mass-init-data}.\\

We perform a perturbation argument. First of all, we linearize around  any equilibrium point $f^0$ of \eqref{stationary-equation}. Let $f^0(p^0)=0$. Assume, without loss of generality, that $p^0=0$. We seek solutions that are perturbations of this equilibrium, as
$$f(x,t)=f^0(x)+g(x,t),$$
where $g$ is a solution of the problem
\begin{equation}\label{eq_g}
g_t=Lg+N(g).
\end{equation}
Here $L$ is the linearized operator (studied in Section \ref{section-linearized}) given by
\begin{equation}\label{op_L}
Lg:=g_{xx}- g_x(p^0)\left[\delta_{x=p^0-a}-\delta_{x=p^0+a}\right] +g(p^0) \left[\delta'_{x=p^0-a}-\delta'_{x=p^0+a}\right],
\end{equation}
and  $N(g)$ contains  the rest of non-linear terms. Note that the free boundary $p(t)$ does not appear in $L$. We first show that, apart from the zero eigenvalue, the rest of eigenvalues of $L$ are strictly negative and isolated, with spectral gap given by $\hat\gamma>0$. Thus the linear problem $g_t=Lg$ has a very nice solution coming from semigroup theory.\\

In order to state our main theorem, we need to fix precisely the function spaces. The basic regularity we need to apply perturbation theory is
$L:X\to X$, with domain $D(L)=Z$, and $N:Y\to X$ for some spaces $Z\subset Y\subset X$. In particular, if we use fractional powers of operators, then we need $Y=X^\alpha$ for some $0<\alpha<1$. These fractional order spaces are simply the $\alpha$-interpolation between $X$ and $Z$.

Naively, we would like to use the spaces
$$Y=H^{r}(-A,B),\quad X=H^{r-2\alpha}(-A,B),\quad Z=H^{r-2\alpha+2}(-A,B),$$
for some $0<r<1/2$. We remark here that the functions we consider will have homogeneous Neumann boundary conditions, so they have a Fourier series expansion $f\sim \sum \hat f_n e_n$. Thus Sobolev spaces $H^r(-A,B)$, $r\in \mathbb R$ can be defined through the norm
$$\norm{f}^2_{H^r}=\sum (1+n^2)^r \hat f_n.$$
The space $X$ is chosen so that $L:X \to X$ is well defined (note that the image of $L$ contains derivatives of $\delta$ functions). In fact for every $0<r<1/2$ one can find $0<\alpha<1$ such that $-3/2>r-2\alpha$.\\

However, in order to get a Lipschitz estimate for the non-linear term of the form
$$\norm{N(g)-N(\tilde g)}_{X}\leq C\norm{g-\tilde g}_Y,$$
we realize that higher regularity is needed near the free boundary $p(t)$. More precisely we will need $g\in \mathcal C^{2,\beta}$. Thus, we will ask the function to (locally) be in some $H^s$ for $3<s<7/2$.
Even in the symmetric case studied by the authors in \cite{price-formation}, higher regularity at the free boundary was required.

To handle this issue,  we introduce a cutoff function $\phi$ around the point $p^0=0$. For the rest of the paper, we fix a small real number $\nu>0$ such that $4\nu<<a$. Assume that $p_I\in(-\nu,\nu)$.
Let $\phi\in\mathcal C^\infty_0(-A,B)$ be a smooth cutoff function such that $\phi\equiv 0$ in $[-A,-2\nu]\cup[2\nu,B]$ and $\phi\equiv 1$ in $[-\nu,\nu]$. Fix $0<r<1/2$, $3<s<7/2$ and consider the space
\be\label{space-X}
X:=\left\{f\in H^{r-2\alpha}(-A,B)\;:\; \phi f\in H^{s-2\alpha}(-A,B)\right\},\ee
with norm
$$\norm{f}_X=\norm{f}_{H^{r-2\alpha}(-A,B)}+\norm{\phi f}_{H^{s-2\alpha}(-A,B)}.$$
On the other hand, we set
$$Z:=\left\{f\in H^{r-2\alpha+2}(-A,B)\;:\; \hat\phi f\in H^{s-2\alpha+2}(-A,B)\right\},$$
with norm
$$\norm{f}_Z=\norm{f}_{H^{r-2\alpha+2}(-A,B)}+\norm{\hat\phi f}_{H^{s-2\alpha+2}(-A,B)}.$$
where $\hat \phi$ is another cutoff function supported on $(-3\nu,3\nu)$ such that $\hat \phi\equiv 1$ on $\mbox{supp } \phi$. Then we have

\begin{lemma}\label{lemma-bounded-operator}
For any $3<s<7/2$, there exist $0<r<1/2$ and $0<\alpha<1$ satisfying all the above requirements such that $L:Z\to X$ is a bounded operator and $Z$ is dense in $X$.
\end{lemma}

The proof of Lemma {lemma-bounded-operator} is contained in \cite{interpolation}, that takes a close look at this kind of functional spaces. However, in order to make the present paper self contained, we rewrite the proof in the Appendix.

For example, take $s=3.4$, $\alpha=0.9$ and $r=0.2$. Thus we can extend $L$ to be an unbounded operator $L:X\to X$ with $D(L)=Z$.
We should ask the initial data to be in the space
$$Y:=X^\alpha=(X,Z)_{\alpha,2},$$
given by the (real) interpolation between $X$ and $Z$. Although we do not have a precise expression for $Y$, in the work \cite{interpolation} we show a characterization of this interpolation space using pseudo-differential operators; in fact, it satisfies  $\bar Y\subset Y\subset \tilde Y$ where
$$\tilde Y:=\left\{f\in H^{r}(-A,B)\;:\; \phi f\in H^{s}(-A,B)\right\},$$
and
$$\bar Y:=\left\{f\in H^{r}(-A,B)\;:\; \bar\phi f\in H^{s}(-A,B)\right\},$$
for some smooth cutoff function  $\bar\phi$ satisfying $\bar\phi\equiv 1$ on the support of $\hat\phi$.
%and $\tilde\phi$ also a smooth cutoff function supported on  $(-\nu,\nu)$ such that $\phi\equiv 1$ on the support of $\tilde\phi$.

%In particular, $s$ is chosen so that if $\phi f\in H^s(-A,B)$, then $f\in \mathcal C^{2,\beta}(-\nu,\nu)$ for some %$0<\beta<1$; as we have mentioned, this embedding is required in order to have a Lipschitz condition for the nonlinear %term $N(g)$ in Section \ref{section-non-linear}.\\

The main result of this work proves the existence of a unique solution $f(x,t)$ for all time $t>0$, if we start with any initial condition $f_I$ {\em{close}} in norm  of the function space $Y$ to a general equilibrium state $f^0$, say
$$\norm{f_I-f^0}_Y\leq \rho,$$
for some $\rho>0$. We also establish exponential decay to a unique stationary state $f_\infty$ (that might not be the same as $f^0$, but it is uniquely determined from the masses $m_1$, $m_2$). Of course, $f_\infty$ will be very near $f^0$, but more significantly, the solution $f(x,t)$ will not leave the neighborhood of size $\rho$, i.e., we have stability; center manifold theory provides a very elegant solution.

\begin{thm}\label{theorem-center-manifold}
Fixed $0<\gamma<\hat \gamma$, for any admissible equilibrium $f^0$ defined in \eqref{equilibrium}, there exist $\rho>0$, $C>0$ such that if we start with initial data
$f_I\in Y$ satisfying
$$\norm{f_I-f^0}_Y\leq \rho,$$
then there exists a unique solution $f(x,t)$ of \eqref{non-linear}-\eqref{mass-init-data} for all time $t>0$ with $f(x,0)=f_I(x)$
 that satisfies
$$f\in\mathcal C^1\lp[0,+\infty):Y\rp,$$
and
\be\label{near-equilibrium}\norm{f(\cdot,t)-f_\infty}_Y\leq C e^{-\gamma t}\norm{f_I-f^0}_Y,\ee
for all $t>0$, where $\hat\gamma$ is the spectral gap given in remark \ref{remark-spectral-gap} and $f_\infty$ is the unique stationary state that satisfies
$$\int_{-A}^{p_\infty} f_\infty =\int_{-A}^{p_I} f_I,\quad \int_{p_\infty}^B f_\infty=\int_{p_I}^B f_I.$$
\end{thm}

An additional interesting ingredient is to see whether the size $\rho$ of the neighborhood can be taken uniformly on $f^0$:

\begin{thm}\label{thm-independent}
The constant $\rho$ in Theorem \ref {theorem-center-manifold} can be taken uniformly, independent of $f^0$, as long as we restrict to the family of equilibria to the set $A_\chi$, where
$$A_\chi:= \left\{f^0 \mbox{admissible}\;:\; \lambda^0\geq \chi\right\},$$
for some $\chi>0$.
\end{thm}

Some comments on the structure of the paper: Section \ref{section-stationary} gives precise formulas for the family of equilibrium points. Section \ref{section-linearized} computes the eigenvalues and eigenfunctions for the linearized operator $L$, while in Section \ref{section-semigroup} we  apply semigroup theory to understand the linearized problem. The main steps in the proof of the theorem are contained in Section \ref{section-non-linear}: first we give the necessary estimates on the non-linear part in order to have local existence; then we show existence and stability of a center manifold so that we have global existence and decay for small solutions, and finally we comment on the stability.

% ==============================================================================================
% ==============================================================================================
\section{The stationary state}\label{section-stationary}

\setcounter{equation}{00}

First we review the construction of the stationary solution from \cite{Lasry-Lions}. Note that, each equilibrium point is determined by two quantities (the masses $m_1$ and $m_2$).

\begin{lemma} \label{lemma-equilibrium}
Given any two positive real numbers $m_1, m_2$, the system (\ref{non-linear})-(\ref{mass-init-data}) has a unique stationary solution $f_\infty(x), p_\infty$ if and only if
\be\label{restriction}
\frac{a}{2A+2B-3a}\leq \frac{m_1}{m_2}\leq \frac{2A+2B-3a}{a}.
\ee
\end{lemma}

\noindent{\emph{Proof: }}
Any equilibrium point of \eqref{non-linear}-\eqref{initial-data} must satisfy
\bee
\left\{\begin{split}
&f^0_{xx}=-\lambda^0\left[ \delta_{x=p^0-a}-\delta_{x=p^0+a}\right]\quad \mbox{in }[-A,B]\times \mathbb R_+\\
&f^0_x(-A)=f^0_x(B)=0\\
& f^0>0\mbox{ if }x\in(-A,p^0), \quad f^0<0\mbox{ if }x\in(p^0,B).
\end{split}\right.
\eee
The solution $f^0$ is a piecewise linear function given by
\bee
f^0_x=\left\{
\begin{split}
0 & \quad\mbox{ if } x\in (0,p^0-a) \\
f^0_x|_{p^0}& \quad\mbox{ if }x\in (p^0-a,p^0+a) \\
0 & \quad \mbox{ if } x\in (p^0+a,A).
\end{split}\right.\eee
The unique stationary state ($f_\infty$, $p_\infty$) is computed with the above formula, if we impose the conservation of mass property,
$$m_1:=\int_{-A}^{p_\infty} f_\infty \,dx,\quad m_2:=-\int_{p_\infty}^{B} f_\infty dx.$$
Since $f_\infty(p_\infty)=0$ and $\lambda_\infty=-\partial_x f_\infty (p_\infty)>0$, then
$$m_1=\lambda_\infty a \; \lp p_\infty-\tfrac{a}{2}+A\rp, \quad m_2=\lambda_\infty a \lp B-p_\infty-\tfrac{a}{2}\rp,$$
thus,
$$\frac{m_1}{m_2}=\frac{p_\infty-\tfrac{a}{2}+A}{B-p_\infty-\tfrac{a}{2}}.$$
Note that the quotient $m_1/m_2$ is an increasing function of $p_\infty$, as expected. From here we get
\be\label{root} p_\infty=\frac{ -a(m_1-m_2)-2Am_2+2Bm_1}{2(m_1+m_2)},\quad \lambda_\infty = \frac{m_1+m_2}{a(-a+A+B)}.\ee
To finish, the condition that $p_\infty\in(-A+a,B-a)$ is equivalent to \eqref{restriction}.
\qed\\

\noindent{\emph{Remark: }} We we say that an equilibrium $f^0$, $f^0(p^0)=0$ as above is \emph{admissible} if it satisfies \eqref{restriction}. Note that the set of admissible equilibria is a smooth family parameterized by $m_1, m_2$, or by $p^0,\lambda^0$, and it has dimension two.
\\
% ===================================================================================================
% ===================================================================================================

\section{Linear stability} \label{section-linearized}

\setcounter{equation}{00}

In this section we look at the eigenvalues of the problem in the interval $[-A,B]$. Let us assume that we are in a situation where an admissible stationary state $f^0$ exists and satisfies $p^0=0$.  Remark that if $m_1>m_2$, then $A>B$, and analogously, if $m_1<m_2$, then $B<A$. On the other hand, $m_1=m_2$, then $A=B$, and the steady state is an odd function, defined on the interval $[-A,A]$.

First, we linearize the equation around the stationary state $f^0,p^0$. If $f$ is a small perturbation of $f^0$, then we can see that it must have a unique root near $p^0$, call it $p(t)$. This perturbation must satisfy
\bee
\left\{\begin{split}
&f_t =f_{xx}-f_x(p(t),t) \left[\delta_{x=p(t)-a}-\delta_{x=p(t)+a}\right]\\
&f(p(t),t)=0\\
&f_x(-A, t)=f_x(B,t)=0.
\end{split}\right.\eee
We write $f=f^0+\epsilon g$, $p=p^0+\epsilon q$, $g=g(x,t)$, $q=q(t)$, differentiate in $\epsilon$ and set $\epsilon=0$. We obtain
\bee
g_t =g_{xx}-\left\{ f_{xx}^0 (p^0)q(t)+g_x(p^0,t) \right\}\left[\delta_{x=p^0-a}-\delta_{x=p^0+a}\right]
 - f^0_x(p^0) \left[\delta'_{x=p^0-a}-\delta'_{x=p^0+a}\right]q(t),
\eee
together with the boundary conditions $g_x  (-A, t)  =g_x(B,t)=0$. Use that $f_{xx}^0 (p^0)=0$, and set $\lambda^0=-f^0_x(p^0)$ so that
\be\label{eq100}
\left\{\begin{split}
&g_t =g_{xx}- g_x(p^0,t)\left[\delta_{x=p^0-a}-\delta_{x=p^0+a}\right] + \lambda^0 \left[\delta'_{x=p^0-a}-\delta'_{x=p^0+a}\right]q(t),\\
&g_x  (-A, t)  =g_x(B,t)=0.
\end{split}\right.\ee
On the other hand, we linearize the condition $f(p(t),t)=0$ and this gives
$$q(t)=\frac{g(p^0)}{\lambda^0},$$
so when we substitute the above formula in \eqref{eq100}, we get the linear equation
\be\label{linearized-equation}
\left\{\begin{split}
&g_t =g_{xx}- g_x(p^0,t)\left[\delta_{x=p^0-a}-\delta_{x=p^0+a}\right] + g(p^0) \left[\delta'_{x=p^0-a}-\delta'_{x=p^0+a}\right],\\
&g_x  (-A, t)  =g_x(B,t)=0.
\end{split}\right.\ee

Note that Lemma \ref{lemma-bounded-operator} assures that the operator $L$ is well defined.\\

In the following, we show that the pointwise spectrum of the operator $L$ consists only of real, non-positive eigenvalues.

\begin{prop}\label{prop-eigenvalues-complex}
Consider the operator
$$Lg:=g_{xx}- g_x(0)\left[\delta_{x=-a}-\delta_{x=a}\right] +g(0) \left[\delta'_{x=-a}-\delta'_{x=a}\right],$$
defined on the space $X$ with $x\in[-A,B]$.
Its pointwise spectrum can be described by
\begin{itemize}
\item Zero is an eigenvalue with eigenspace of dimension two. Two linearly independent eigenfunctions are given by
\bee g_0(x)=\left\{\begin{split}
&x, & x\in(-a,a) \\
&a, &x \in(a,B)\\
&-a, &\quad  x\in (-A,-a)
\end{split}\right.,\quad
h_0(x)=\left\{\begin{split}
&1, & x\in(-a,a) \\
&2, &x \in (a,B)\\
&2, &\quad  x\in (-A,-a).
\end{split}\right.
\eee

\item The rest of the eigenvalues must be of the form $$\mu=-\alpha^2,\quad \mbox{for}\quad \alpha=\frac{n\pi}{\alpha},\frac{n\pi}{2A-a}, \frac{n\pi}{2B-a}, \quad \mbox{for some } n\in \mathbb Z.$$
Every eigenspace has finite dimension.
\end{itemize}
\end{prop}

\noindent\emph{Proof: } Ideas come from the previous work of the authors \cite{price-formation} on the symmetric case.
Let  $\mu\in\mathbb C$, and $g\in X$ such that $Lg=\mu g$, with boundary conditions
\be\label{Neumann-complex}g_x(-A)  =g_x(B)=0.\ee
More precisely,
$$g_{xx}- g_x(0)\left[\delta_{x=-a}-\delta_{x=a}\right] +g(0) \left[\delta'_{x=-a}-\delta'_{x=a}\right] =\mu g.$$
We see that the function $g$ must be discontinuous at $-a$ and $+a$ with jump,
\be\label{matching1}\begin{split}
&g(a^+)-g(a^-)=g(0),\\
&g(-a^+)-g(-a^-)=-g(0),
\end{split}\ee
and that its derivative too, with jump
\be\label{matching2}
\begin{split}
&g'(a^+)-g'(a^-)=-g'(0),\\
&g'(-a^+)-g'(-a^-)=g'(0).
\end{split}\ee
It is clear that, except at the points $x=\pm a$, the function $g$ must satisfy $g_{xx}=\mu g$. Thus, our main trick is to find our eigenfunctions in a similar way as one finds eigenvalues and eigenfunctions for the Laplacian using Fourier series.

Let $\alpha=+\sqrt{\mu}\in\mathbb C$, $\alpha\neq 0$.  We set
\bee\left\{
\begin{split}
&g_1(x)=c_1 e^{\alpha x}+c_2 e^{-\alpha x} &\quad\mbox{ in }(-a,a),\\
&g_2(x)=d_1 e^{\alpha x}+d_2 e^{-\alpha x}&\quad\mbox{ in }(a,B),\\
&g_3(x)=e_1 e^{\alpha x}+e_2e^{-\alpha x}&\quad\mbox{ in }(-A,-a).\\
\end{split}\right.\eee
This $g$ must still satisfy the boundary conditions \eqref{Neumann-complex} and the matching conditions \eqref{matching1}-\eqref{matching2} as given above. When we impose \eqref{Neumann-complex}, it is easy to see that $g$ can be written simply as
\bee\left\{
\begin{split}
&g_1(x)=c_1 e^{\alpha x}+c_2 e^{-\alpha x} &\quad\mbox{ in }(-a,a),\\
&g_2(x)=d \left[e^{-2\alpha B+\alpha x}+ e^{-\alpha x}\right]&\quad\mbox{ in }(a,B),\\
&g_3(x)=e \left[e^{2\alpha A+\alpha x}+e^{-\alpha x}\right]&\quad\mbox{ in }(-A,-a).\\
\end{split}\right.\eee
The matching conditions at $x=a$ imply that
\bee\left\{\begin{split}
&d\left[ e^{2\alpha B+\alpha a} + e^{ -\alpha a} \right] +c_1\left[ -e^{\alpha a}-1\right]+c_2\left[-e^{-\alpha a}-1\right]=0,\\
&d\left[ e^{-2\alpha B+\alpha a} - e^{ -\alpha a} \right]+c_1\left[ -e^{\alpha a}+1\right]+c_2\left[e^{-\alpha a}-1\right]=0.
\end{split}\right.\eee

When we add and substract the two equations above we get:
\be\label{1-2}\begin{split}
& c_1 e^{\alpha a} +c_2 +d  e^{-2\alpha B+\alpha a}=0,\\
& c_1 +  c_2 e^{-\alpha a}  -d  e^{-\alpha a}=0.
\end{split}\ee

On the other hand, we can substitute the matching conditions at $x=-a$ to obtain:
\bee\begin{split}
&c_1 \left[e^{-\alpha a}+ 1\right]+c_2 \left[e^{\alpha a}+1\right]-e\left[ e^{2\alpha A-\alpha a} + e^{ \alpha a} \right]=0,\\
&c_1\left[ e^{-\alpha a}-1\right]+c_2\left[-e^{\alpha a}+1\right]-e\left[ e^{2\alpha A-\alpha a} - e^{\alpha a} \right]=0.
\end{split}\eee
Then, adding and substracting the previous two equations:
\be\label{3-4}\begin{split}
&e^{-\alpha a}c_1+c_2-e e^{2\alpha A -\alpha a}=0, \\
&c_1 +c_2 e^{\alpha a}-e e^{\alpha a}=0.
\end{split}\ee

Then we see that \eqref{1-2} together with \eqref{3-4} become an homogeneous linear system of four equations with four unknowns. Performing row reduction in the system, we can compute the determinant of the coefficient matrix, and indeed, it is a (non-zero) multiple of
$$\lp e^{-\alpha a}-e^{-2\alpha B}\rp \lp 1-e^{-2\alpha a}\rp \lp e^{-2\alpha A }-e^{\alpha a}\rp.$$
It is clear that system \eqref{1-2}-\eqref{3-4} only has  the trivial solution unless one of the three factors above vanishes. In this case, we must have that the real part of $\alpha$ is zero, and the imaginary part must take the values:
$$\mathcal{I}(\alpha)=\frac{n\pi}{\alpha},\frac{n\pi}{2A-a}, \frac{n\pi}{2B-a}, \quad n\in \mathbb Z.$$
The eigenvalues must necessarily be of the form $\mu=\alpha^2$, which are all real and negative. The eigenspace for each eigenvalue is finite dimensional and depends on the number of solutions of the linear system above.

To complete the proof of the proposition, we need to check the zero eigenspace. Then we seek functions $g$ such that
$$g''- g'(0)\left[\delta_{x=-a}-\delta_{x=a}\right] +g(0) \left[\delta'_{x=-a}-\delta'_{x=a}\right] =0.$$
In particular only the following functions $g_0$ and $h_0$ (and a linear combination of those) will be suitable solutions:
\bee
\left\{\begin{split}
& g_0(0) =0 & \textrm{and (\ref{matching2}) is satisfied, }  \\
& h_0'(0) =0 & \textrm{and (\ref{matching1}) is satisfied}.
\end{split}\right.
\eee
The corresponding zero eigenspace is given by a linear combination of $g_0$ and $h_0$:
 \bee g_0(x)=\left\{\begin{split}
&x, & x\in(-a,a) \\
&a, &x \in(a,B)\\
&-a, &\quad  x\in (-A,-a)
\end{split}\right.,\quad
h_0(x)=\left\{\begin{split}
&1, & x\in(-a,a) \\
&2, &x \in (a,B)\\
&2, &\quad  x\in (-A,-a)
\end{split}\right.
.\eee
It is an easy computation to check that $g_0$ and $h_0$ are linearly independent.

\qed

It is clear from the method above how to compute the eigenfunctions. For completeness here, we show the explicit calculations in the case $A=B=1$.

\begin{prop}\label{prop-eigenfunctions}
Consider the operator
$$Lg:=g_{xx}- g_x(0)\left[\delta_{x=-a}-\delta_{x=a}\right] +g(0) \left[\delta'_{x=-a}-\delta'_{x=a}\right],$$
defined on the space $X$ with $x\in[-1,1]$.
All the non-zero  eigenvalues of $L$ are given by:
\begin{enumerate}

\item For all $n$ such that
$$\frac{2n}{a}- \frac{1}{2} \in\mathbb Z,$$
we have
$$\mu_n=-(\alpha_n)^2\quad\mbox{for}\quad \alpha_n:=\frac{2n\pi}{a}$$
with corresponding one dimensional eigenspace with eigenfunction given by
\bee g_n(x)=\left\{\begin{split}
&\sin(\alpha x)  &\quad x\in(-a,a) \\
&0&\quad x \in(a,1)\\
&0, &\quad  x\in (-1,-a).
\end{split}\right.
\eee

\item For all $n$ such that
$$
\frac{(2n+1)}{a}- \frac{1}{2}\in\mathbb Z,$$
we have
$$
\mu_n=-(\alpha_n)^2\quad\mbox{for}\quad \alpha_n:=\frac{(2n+1)\pi}{a}
$$
with corresponding two dimensional eigenspace given by the linear combination of the following functions

\bee g_n(x)=\left\{\begin{split}
&\sin(\alpha x)  &\quad x\in(-a,a) \\
&2 \sin(\alpha x) &\quad x \in(a,1)\\
&2 \sin(\alpha x) , &\quad  x\in (-1,-a).
\end{split}\right.
\eee

and

\bee h_n(x)=\left\{\begin{split}
&\cos(\alpha x)  &\quad x\in(-a,a) \\
&0 &\quad x \in(a,1)\\
&0, &\quad  x\in (-1,-a).
\end{split}\right.
\eee

\item For all $n$ such that
$$
\frac{(2n+1)}{a} - \frac{1}{2} \not\in\mathbb Z,$$
we have
$$
\mu_n=-(\alpha_n)^2\quad\mbox{for}\quad \alpha_n:=\frac{(2n+1)\pi}{a}
$$
with corresponding one dimensional eigenspace generated by

\bee g_n(x)=\left\{\begin{split}
&\cos(\alpha x)  &\quad x\in(-a,a) \\
&0 &\quad x \in(a,1)\\
&0, &\quad  x\in (-1,-a).
\end{split}\right.
\eee

\item For all $n$ such that
$$
\frac{2n}{a} - \frac{1}{2} \not\in \mathbb Z, \quad \textrm{and} \quad \frac{2n}{a} \not\in \mathbb Z,$$
we have
$$
\mu_n=-(\alpha_n)^2\quad\mbox{for}\quad \alpha_n:=\frac{2n\pi}{a}
$$
with corresponding one dimensional eigenspace generated by

\bee g_n(x)=\left\{\begin{split}
&\sin(\alpha x)  &\quad x\in(-a,a) \\
&0 &\quad x \in(a,1)\\
&0, &\quad  x\in (-1,-a).
\end{split}\right.
\eee

\item For all $n$ such that
$$
\frac{2n}{a} - \frac{1}{2} \not\in\mathbb Z , \quad \textrm{and} \quad \frac{2n}{a} \in\mathbb Z,$$
we have
$$
\mu_n=-(\alpha_n)^2\quad\mbox{for}\quad \alpha_n:=\frac{2n\pi}{a}
$$
with corresponding two dimensional eigenspace generated by

\bee g_n(x)=\left\{\begin{split}
&\sin(\alpha x)  &\quad x\in(-a,a) \\
&0 &\quad x \in(a,1)\\
&0, &\quad  x\in (-1,-a),
\end{split}\right.
\eee

and

\bee h_n(x)=\left\{\begin{split}
&\cos(\alpha x)  &\quad x\in(-a,a) \\
& 2 \cos(\alpha x) &\quad x \in(a,1)\\
& 2\cos(\alpha x), &\quad  x\in (-1,-a).
\end{split}\right.
\eee

\item For all $n$ such that
$$
\frac{2n}{2-a} - \frac{1}{2} \not\in\mathbb Z, \quad \textrm{and} \quad \frac{2n}{2-a}a \not\in\mathbb Z,$$
we have
$$
\mu_n=-(\alpha_n)^2\quad\mbox{for}\quad \alpha_n:=\frac{2n \pi}{2-a}
$$
with corresponding two dimensional eigenspace generated by
\bee g_n(x)=\left\{\begin{split}
&-\frac{\sin(\alpha(1-a))}{2(1-\cos(\alpha a))}\sin(\alpha x)+\frac{\cos(\alpha(1-a))}{2(1+\cos(\alpha a))}\cos(\alpha x),  &\quad x\in(-a,a) \\
&  \cos(\alpha (1-x)), &\quad x \in(a,1)\\
& 0, &\quad  x\in (-1,-a),
\end{split}\right.
\eee
and
\bee h_n(x)=\left\{\begin{split}
&\frac{\sin(\alpha(1-a))}{2(1-\cos(\alpha a))}\sin(\alpha x)+\frac{\cos(\alpha(1-a))}{2(1+\cos(\alpha a))}\cos(\alpha x),  &\quad x\in(-a,a) \\
&  0, &\quad x \in(a,1)\\
& \cos(\alpha (1+x)), &\quad  x\in (-1,-a).
\end{split}\right.
\eee
\end{enumerate}
Moreover, for all $n\geq 1$, the eigenfunctions have zero mass, i.e.,
\be\label{mass-eigenfunctions}\int_{-1}^0 g_n\;dx=\int_{-1}^0 h_n\;dx =0,\quad \int_0^1 g_n\;dx= \int_0^1 h_n \;dx=0.\ee

\end{prop}

\noindent{\emph{Proof: }}
We have seen that all the eigenvalues are of the form $\mu=-\alpha^2\leq 0$.
The construction of the eigenfunction $g$ is done piecewise in each of these three intervals: $(-a,a)$, $(-1,-a)$, and $(a,1)$, using Fourier series. We seek
\bee\left\{
\begin{split}
&g_1(x)=c_1\sin(\alpha x)+c_2\cos(\alpha x) &\quad\mbox{ in }(-a,a),\\
&g_2(x)=d_1\sin(\alpha x)+d_2\cos(\alpha x)&\quad\mbox{ in }(a,1),\\
&g_3(x)=e_1\sin(\alpha x)+e_2\cos(\alpha x)&\quad\mbox{ in }(-1,-a),\\
\end{split}\right.\eee
with boundary conditions
\be\label{zero-neumann}g'_3(-1)=g'_2(1)=0.\ee

The matching conditions \eqref{matching1}, \eqref{matching2} are rewritten as
\be\begin{split} \label{jumps}
&g_2(a)-g_1(a)=g_1(0),\\
&g_1(-a)-g_3(-a)=-g_1(0),\\
&g_2'(a)-g'_1(a)=-g'_1(0),\\
&g'_1(-a)-g'_3(-a)=g'_1(0).
\end{split}\ee
The zero Neumann boundary conditions \eqref{zero-neumann} give that
\bee\begin{split}
d_1 \cos(\alpha)-d_2\sin(\alpha)=0,\\
e_1 \cos(\alpha)+e_2\sin(\alpha)=0.
\end{split}\eee

Consider first the case $\cos(\alpha) =0$.

This implies $d_2 = e_2 =0$ and from the matching conditions (\ref{jumps}) we get the system
\be\left\{
\begin{split}\label{case00}
d_1 \sin(\alpha a ) - c_1 \sin (\alpha a ) - c_2 \cos( \alpha a) = c_2, \\
d_1 \cos(\alpha a ) - c_1 \cos (\alpha a ) + c_2 \sin( \alpha a) = - c_1, \\
-c_1 \sin(\alpha a ) + c_2 \cos (\alpha a ) + e_1 \sin( \alpha a) =  - c_2, \\
c_1 \cos (\alpha a ) + c_2 \sin ( \alpha a)  -e_1 \cos ( \alpha a) = c_1.
\end{split} \right.\ee

Consider now the case $\sin(\alpha a ) \neq 0$: eliminating the constant $c_1$ from the first and second equation in (\ref{case00}), we get
$$d_1 \sin(\alpha a ) =0,$$
which is satisfied only if $d_1=0$. With a similar computations in the third and fourth equations, we get $e_1 \sin(\alpha a ) =0$ which is satisfied only for $e_1=0$. Therefore system (\ref{case00}) reduces to
\bee\left\{
\begin{split}
 - c_1 \sin (\alpha a ) - c_2 \cos( \alpha a) = c_2, \\
 - c_1 \cos (\alpha a ) + c_2 \sin( \alpha a) = - c_1, \\
-c_1 \sin(\alpha a ) + c_2 \cos (\alpha a )  =  - c_2, \\
c_1 \cos (\alpha a ) + c_2 \sin ( \alpha a) = c_1,
\end{split} \right.\eee
which is the unique solution $c_1 = c_2=0$.

Therefore consider now $\alpha$ such that $   \sin(\alpha a ) = 0$: this yields to
\bee\left\{
\begin{split}
 - c_2 \cos( \alpha a) = c_2, \\
d_1 \cos(\alpha a) - c_1 \cos (\alpha a ) = - c_1, \\
c_2 \cos (\alpha a )  =  - c_2, \\
c_1 \cos (\alpha a ) - e_1 \cos ( \alpha a) = c_1,
\end{split} \right.\eee
which has solutions:
$$
d_1 = e_1 = c_2 =0, \quad \textrm{for} \quad  \cos (\alpha a ) = 1,$$
and
$$
 d_1 =  2c_1, \; e_1 = 2 c_1 , \quad \textrm{for} \quad  \cos (\alpha a ) = -1.$$
This proves Assertion (1) and (2).

Consider now the case $\cos(\alpha)  \neq 0$.

We can write
$$
d_1 = d_2 \frac{ \sin( \alpha)}{ \cos (\alpha) }, \quad e_1 = -e_2 \frac{ \sin( \alpha)}{ \cos (\alpha) },
$$
which implies that from the matching conditions (\ref{jumps}) we get the following system
\be\left\{
\begin{split}\label{case11}
d_2 \cos(\alpha (1-a) ) - c_1 \sin (\alpha a ) - c_2 \cos( \alpha a) = c_2, \\
d_2 \sin(\alpha (1-a)  ) - c_1 \cos (\alpha a ) + c_2 \sin( \alpha a) = - c_1, \\
-c_1 \sin(\alpha a ) + c_2 \cos (\alpha a ) - e_2 \cos( \alpha (1-a)) =  - c_2, \\
c_1 \cos (\alpha a ) + c_2 \sin ( \alpha a)  + e_2 \sin ( \alpha (1-a)) = c_1.
\end{split} \right.\ee

Suppose  $   \sin(\alpha a ) = 0$ (we remind the reader that we are under the assumption $\cos(\alpha)  \neq 0$): the above system reduces to

\be\left\{
\begin{split}  \label{case110}
d_2 \cos(\alpha a ) \cos(\alpha) - c_2 \cos( \alpha a) = c_2, \\
d_2 \sin(\alpha ) \cos(\alpha a )  - c_1 \cos (\alpha a )  = - c_1, \\
 c_2 \cos (\alpha a ) - e_2 \cos( \alpha)\cos(\alpha a )  =  - c_2, \\
c_1 \cos (\alpha a ) + e_2 \sin ( \alpha) \cos(\alpha a )  = c_1.
\end{split} \right.\ee

If $\cos(\alpha a ) = -1$,  system (\ref{case110}) has a unique solution
$$
d_2 = e_2 = c_1 =0,
$$
and Assertion (3) is proven. \\
Otherwise, if $\cos(\alpha a ) = 1$, we separate the case when $\sin(\alpha) =0$ and $\sin(\alpha)  \neq 0$. \\
For $\sin(\alpha) \neq 0$ system (\ref{case110}) has solution $d_2 =  e_2 = c_2 =0$, which proves Assertion (4).

For $\sin(\alpha) = 0$ system (\ref{case110}) has a unique solution
$$
d_2 = 2 c_2, e_2 = 2 c_2,
$$
and Assertion (5) follows.

Suppose  $\sin(\alpha a ) \neq  0$ and again $\cos(\alpha)  \neq 0$. Eliminating $c_1$ from the first and second equation of \eqref{case11} we get
$$d_2 \cos (\alpha(1-a))=d_2\cos(\alpha).$$
On the other hand, eliminating $c_2$ gives
$d_2\sin(\alpha)=-d_2 \sin(\alpha(1-a))$, which implies, that $d_2=0$ or $\alpha= \frac{2n\pi}{2-a}$. Similar computations using the third and fourth equation of \eqref{case11} yield to $e_2=0$ or $\alpha=\frac{2n\pi}{2-a}$. However, we can check that $d_2=e_2=0$ does not produce any non-trivial solution. For $\alpha=\frac{2n\pi}{2-a}$, combining the first and third equations of \eqref{case11} we get that
$$d_2 \cos (\alpha (1-a) ) + e_2 \cos( \alpha (1-a)) = 2 c_2 ( 1 + \cos(\alpha a)),$$
which implies
$$c_2 =  \frac{1}{2}  \cos (\alpha (1-a) ) \frac{ d_2 + e_2}{ ( 1 + \cos(\alpha a))}.$$
Similarly, from second and fourth equation of (\ref{case11}) we get
$$c_1 =  \frac{1}{2}  \sin (\alpha (1-a) ) \frac{ e_2 - d_2}{ ( 1 + \cos(\alpha a))}.$$
Assertion (6) is proven.

Summarizing, the structure of the eigenvalues follows the following scheme:
\bee\left\{
\begin{split}
 & \cos(\alpha)=0    \left\{  \begin{split}
 & \sin(\alpha a ) \neq  0 \quad \textrm{no solutions} \\
 & \sin(\alpha a ) =  0  \left\{ \begin{split} &\cos(\alpha a)=1   \quad \quad\textrm{(Assertion 1)} \\
& \cos(\alpha a)=- 1 \quad \quad  \textrm{(Assertion 2)},  \end{split}\right.   \end{split}\right. \\
  & \cos(\alpha)\neq 0
  \left\{  \begin{split}  &\sin(\alpha a ) =  0   \left\{  \begin{split}  &\cos(\alpha a ) =  -1 \quad \textrm{(Assertion 3)} \\
 & \cos(\alpha a ) =  1   \left\{ \begin{split} &\sin(\alpha)=0 \quad \quad \textrm{(Assertion 5)} \\
 &\sin(\alpha)\neq0 \quad \quad \textrm{(Assertion 4)},  \end{split}\right.
 \end{split} \right. \\
 & \sin(\alpha a ) \neq 0   \quad \textrm{(Assertion 6).}
 \end{split} \right.
 \end{split} \right.\eee

\bigskip

\begin{rem}\label{remark-spectral-gap}
We write the spectral gap as
\bee
\hat \gamma:=  \min\left\{ \lp\frac{2\pi}{2A-a}\rp^2, \lp\frac{2\pi}{2B-a}\rp^2, \lp\frac{\pi}{a}\rp^2\right\}.
\eee
\end{rem}

% ==============================================================================================================
% ==============================================================================================================

\section{Semigroup theory} \label{section-semigroup}

\setcounter{equation}{00}
In the following, we study the properties of the linear part of our equation \eqref{eq_g}, given by $g_t = Lg,$
where the operator $L$ is defined as in \eqref{op_L}. In particular we will show that it generates an analytic semigroup, together with time dependent decay estimates. We close the section with a very explicit characterization of the $ker\; (L)$, where
$$
\ker(L) = \textrm{span}\langle g_0, h_0\rangle,
$$
with $g_0$ and $h_0$ defined as in  Proposition {\ref{prop-eigenfunctions}}.\\

To begin with, we review some important concepts from semigroup theory. Standard references are the first chapter of the book \cite{Henry:parabolic} and also \cite{Pazy:book}. Let $X$ be a Banach space with norm $\norm{\cdot}_X$ and  let $L$ be a linear operator on $X$ that has domain and range in $X$. Denote by $\sigma(L)\subset \mathbb C$ its spectrum.

We say that $L$ is a \emph{sectorial operator} if it is a closed densely defined operator such that, for some $\phi\in(0,\pi/2)$ and some $M\geq 1$ and real $a$, the sector
$$S_{a,\phi}:=\left\{\lambda\in\mathbb C : \phi\leq \abs{\arg (\lambda-a)}\leq \pi,\lambda\neq a\right\}$$
is in the resolvent set of $L$ and
$$\norm{(\lambda I-L)^{-1}}_X\leq M/\abs{\lambda-a}\quad \mbox{for all }\lambda\in S_{a,\phi}.$$
An \emph{analytic semigroup} on a Banach space $X$ is a family of continuous linear operators on $X$, $\{T(t)\}_{t\geq 0}$, satisfying
\begin{enumerate}
\item $T(0)=I$, $T(t)T(s)=T(t+s)$ for $t\geq 0$, $s\geq 0$.
\item $T(t)x\to x$ as $t\to 0^+$, for each $x\in X$.
\item $t \to T(t)x$ is real analytic on $0<t<\infty$ for each $x\in X$.
\end{enumerate}
The \emph{infinitesimal generator} $L$ of this semigroup is defined by
$$Lx=\lim_{t\to 0^+} \frac{T(t)x-x}{t},$$
its domain $D(L)$ consisting of all $x\in X$ for which this limit exists. And viceversa, it is well known that
(see Theorem 1.3.4 in \cite{Henry:parabolic}, for instance) if $L$ is a sectorial operator, then $-L$ is the infinitesimal generator of an analytic semigroup $\{T(t)=e^{-tL}\}_{t\geq 0}$; this semigroup gives the solution of the ODE $g_t+Lg=0$. Moreover, if $\Re \sigma(L)>\gamma$, then for $t>0$, we have the bounds
$$\norm{e^{-Lt}}_X\leq C e^{-\gamma t},\quad \norm{Le^{-Lt}}_X\leq \frac{C}{t}e^{-\gamma t},$$
for some constant $C$. This implies in particular that $\textrm{Range}\;(e^{-Lt}) \subset D(L)$.\\

Suppose that $L$ is a sectorial operator and $\Re\sigma(L)>0$. Then for any $\alpha >0$ we define the $(-\alpha)$-\emph{fractional power} as
$$L^{-\alpha}=\frac{1}{\Gamma(\alpha)}\int_0^\infty t^{\alpha-1}e^{-Lt}\,dt.$$
The operator $L^{-\alpha}$ is a bounded linear operator on $X$ which is one-to-one.
For $\alpha\geq 0$, $L^\alpha$ is defined to be the inverse of $L^{-\alpha}$, $D(L^\alpha)=R(L^{-\alpha})$.   $L^\alpha$ is a closed, densely defined operator. We define
$$X^\alpha:=D(L^\alpha),$$
with the graph norm
$$\norm{x}_{X^\alpha}:=\norm{L^\alpha x}_X.$$
It is seen in the last section that the space $X^\alpha$ is precisely the $\alpha$-interpolation between the spaces $X$ and $D(L)$.\\

If the spectrum of the operator $L$ has a good structure, then we can split the space $X$.  A set $\sigma\subset \sigma(L)\cup \{\infty\}=:\hat \sigma(L)$ is a spectral set if both $\sigma$ and $\hat\sigma (L)\backslash \sigma$ are closed in the extended plane $\mathbb C\cup \{\infty\}$.

\begin{thm} [Theorem 1.5.2. in \cite{Henry:parabolic}]\label{thm-splitting}
Suppose $L$ is a closed linear operator in $X$ and suppose $\sigma_1$ is a bounded spectral set, and $\sigma_2=\sigma(L)\backslash \sigma_1$ so $\sigma_2\cup\{\infty\}$ is another spectral set. Let $E_1$, $E_2$ be the projections associated with these spectral sets, and $X_j=E_j(X)$, j=1,2. Then $X=X_1\oplus X_2$, the $X_j$ are invariant under $L$, and if $L_j$ is the restriction of $L$ to $X_j$, then
$$L_1:X_1\to X_1\mbox{ is bounded}, \sigma(L_1)=\sigma_1,$$
$$D(L_2)=D(L)\cap X_2\mbox{ and }\sigma(L_2)=\sigma_2.$$
\end{thm}

The following lemmas will allow us to compare our operator $L$ to the standard Laplacian. The first one will be used in order to prove that our linear operator $L$ is sectorial, while the second gives that the spectrum of $L:=\Delta+B$ consists only of eigenvalues.

\begin{lemma}[Theorem 1.4.5 in \cite{Henry:parabolic}]\label{compare-sectorial}
If $M$ is a sectorial operator with $\Re{\sigma(M)}>0$ and if $B$ is a linear operator such that $BM^{-\alpha}$ is bounded on $X$ for some $0<\alpha<1$, then $M+B$ is sectorial.
\end{lemma}

\begin{lemma}[Chapter II, section 5.13, and Chapter V, Corollary 1.15 in \cite{Engel}]\label{Engel-lemma}
Let $M$ be an operator having compact resolvent and let $B$ be a bounded operator on $X$ such that the resolvent set of $M+B$ is non-empty. Then $M+B$ has compact resolvent. As a consequence, the spectrum of $M+B$ consists only of eigenvalues.
\end{lemma}

The next lemma is used in order to define the domain of a fractional power of an operator $X^\alpha(=Y)$, and to get a time decay estimate for the semigroup. This space is important because it is the space of admissible initial conditions; the lemma tells us that it is enough to consider fractional powers of the Laplacian to define it. However, a more explicit characterization of $Y$ in terms of usual norms is desirable - this is done in \cite{interpolation}.

\begin{lemma}[Theorem 1.5.4. in \cite{Henry:parabolic}]\label{lemma-estimates}
Suppose $L$ is a sectorial operator, $\sigma_1$ a bounded spectral set for $L$, $\sigma_2=\sigma(L)\backslash \sigma_1$, $\Re\sigma_2>\gamma$ and $X=X_1\oplus X_2$ is the corresponding decomposition.
Assume also that $M$ is a sectorial operator with $D(M)=D(L)$, $\Re\sigma (M)>0$, $(M-L)M^{-\alpha}$ is bounded for some $0<\alpha<1$. Then using the norm $\norm{x}_{X^\alpha}:=\norm{M^\alpha x}$, $0\leq \alpha\leq 1$, for $x\in X_2\cap D(M^{\alpha})$ and $t>0$,
$$\norm{e^{-L_2 t}x}_{X^\alpha}\leq C_1\norm{x}_Xt^{-\alpha}e^{-\gamma t},$$
$$\norm{e^{-L_2 t}x}_{X^\alpha}\leq C_1\norm{x}_{X^\alpha} e^{-\gamma t},$$
for some positive constant $C_1$.
\end{lemma}

%\bigskip

%Now we look at our operator $L$ defined as in \eqref{op_L}. Given $5/2<s<3$, $0<r<1/2$, we can find $0<\alpha<1$ such %that
%$-2<r-2\alpha<-3/2$. As we have mentioned, we define the spaces $Z\subset Y\subset X$ as
%\begin{align}
%X&=\left\{f\in H^{r-2\alpha}(-A,B)\;:\; \phi f\in H^{s-2\alpha}(-A,B)\right\},\\
%Z&=\left\{f\in H^{r-2\alpha+2}(-A,B)\;:\; \phi f\in H^{s-2\alpha+2}(-A,B)\right\}.
%\end{align}
%Then $L:X \to X$ and $D(L)=Z$. Let $Y:=X^\alpha$.

Consider now our linear problem
\be\label{problem-linear}\left\{\begin{split}
&g_t = Lg,\\
&g(x,0)=g_I,\\
&g_x(-A,t)=g_x(B,t)=0.
\end{split}\right.\ee
where $L:X\to X$ is defined as in \eqref{op_L}, and the initial condition $g_I\in Y$.
Using the results recalled above, we are able to prove now the next proposition:

\begin{prop}\label{prop-basic}
Let $\sigma_1=\{0\}$ and  $\sigma_2$ be the rest of non-zero eigenvalues of $L$, defined as in \eqref{op_L}. Let $E_1$, $E_2$ be the projections associated with these spectral sets, and $X_j=E_j(X)$, j=1,2. Then we have the splitting $X=X_1\oplus X_2$, where the $X_j$ are invariant under $L$, and if $L_j$ is the restriction of $L$ to $X_j$, then
$$L_1:X_1\to X_1\mbox{ is bounded }, \sigma(L_1)=\sigma_1,$$
$$D(L_2)=D(L)\cap X_2\mbox{ and }\sigma(L_2)=\sigma_2.$$
In our case, $E_1=\ker (L)$, $L_1\equiv 0$. Moreover, $-L_2$ is sectorial  generates an analytic semigroup $T(t)=e^{L_2t}$.
This semigroup satisfies the estimates
\be\label{estimates-semigroup}\begin{split}
&\norm{T(t)x}_Y\leq C_1\norm{x}_X t^{-\alpha}e^{-\gamma t}\\
&\norm{T(t)x}_Y\leq C_1\norm{x}_Y e^{-\gamma t}.
\end{split}\end{equation}
for any $0<\gamma<\hat\gamma$.
\end{prop}

\noindent{\emph{Proof: }}
Proposition \ref{prop-eigenfunctions} computes the spectrum of $L$. Apart from the zero eigenvalue (with eigenspace of dimension two), the operator $-L$ has a discrete and positive spectrum. Thus Theorem \ref{thm-splitting} splits the space $X$: let $\sigma_1=\{0\}$ and $\sigma_2$ be the rest of eigenvalues. Then  $X=X_1\oplus X_2$
where $X_1=\ker(L)$, that satisfies $\dim(X_1)=2$.

Next, we write $-L_2=-\Delta+B_{|_{X_2}}$ where
$$B=\left.\lp g_x(p^0)\left[\delta_{x=p^0-a}-\delta_{x=p^0+a}\right]-g(p^0)\left[\delta'_{x=p^0-a}-\delta'_{x=p^0+a}\right]\rp\right. .$$
However, $-\Delta_{|_{X_2}}$ contains only positive eigenvalues. Thus Lemma \ref{compare-sectorial} applied to $M=-\Delta_{|_{X_2}}$ and $B$ as above gives that $-L_2$ is sectorial, once we note that $D(-\Delta_{|_{X_2}}^\alpha) \subset D(B)$. In fact, since $D(-\Delta_{|_{X_2}}^\alpha) \subset D(B)$, closed graph theorem implies that $B\Delta_{|_{X_2}}^{-\alpha}$ is bounded. \\

Next, we use Lemma \ref{lemma-estimates} for $M=-\Delta_{|_{X_2}}$, and $\beta=\alpha\in(0,1)$, so
$$\norm{e^{L_2 t}x}_Y\leq C_1\norm{x}_Xt^{-\alpha}e^{-\gamma t},$$
$$\norm{e^{L_2 t}x}_Y\leq C_1\norm{x}_Y e^{-\gamma t},$$
for some $C_1$ positive constant.\\
\qed\\

\bigskip
The solution of equation \eqref{problem-linear} can be written now as
$$
g(x,t) = c_I g_0 + d_I h_0 + e^{L_2t}(g_I),
$$
with $c_I  g_0 + d_I h_0 \in X_1$ for some $c_I, d_I\in\mathbb R$ that depend only on the initial condition $g_I$, and $e^{L_2t}(g_I) \in X_2$.\\

Although it is not usually known, in our situation we are able to give a more explicit characterization of $X_1=\ker L$. It will be needed later.

\begin{lemma}\label{characterize-kernel}
Given any $g\in X$, $g=g_1+g_2$, $g_i\in X_i$, $i=1,2$, then in the basis of eigenfunctions $\{g_0,h_0\}$ given in Proposition \ref{prop-eigenfunctions} we can write
$$g_1=c g_0+ d h_0\in X_1,$$
for
\be\label{coefficients}
\begin{split}
&c=\frac{1}{a(a-2A)(a-2B)}\left\{(-a+2A)I_2[g]-(-a+2B)I_1[g]\right\},\\
&d=\frac{1}{a(a-2A)(a-2B)}\left\{ (-\tfrac{a^2}{2}+a B)I_1[g]-(\tfrac{a^2}{2}-aA)I_2[g]\right\},
\end{split}\ee
where $I_1[g]$ and $I_2[g]$ are defined as
$$I_1[g]:=\int_{-A}^0 g\;dx,\quad I_2[g]:=\int_0^B g\;dx.$$
\end{lemma}

\noindent{\emph{Proof: }}
We know that $g=c g_0+d h_0+g_2$, for $g_2\in X_2$ for some $c,d$. If we integrate this expression we obtain
\bee\begin{split}
\int_0^B g\;dx= c\int_0^B g_0\;dx+d\int_0^B h_0\;dx +\int_0^B g_2\;dx,\\
\int_{-A}^0 g\;dx= c\int_{-A}^0 g_0\;dx+d\int_{-A}^0 h_0\;dx +\int_{-A}^0 g_2\;dx.
\end{split}\eee
However, because $g_2\in X_2$ (see \eqref{mass-eigenfunctions})
we know that
$$\int_{-A}^0 g_2\;dx=\int_0^B g_2\;dx=0,$$
thus we conclude from the previous calculation that $c$, $d$ must solve
\bee\begin{split}
I_2[g]= c\int_0^B g_0\;dx+d\int_0^B h_0\;dx, \\
I_1[g]= c\int_{-A}^0 g_0\;dx+d\int_{-A}^0 h_0\;dx.
\end{split}\eee
The integrals of $g_0$, $h_0$ can be explicitly computed; we obtain a system of two equations and two unknowns $c,d$
\bee\begin{split}
I_2[g]= c\lp-\tfrac{a^2}{2}+aB\rp +d(2B-a), \\
I_1[g]= c\lp \tfrac{a^2}{2}-aA\rp+d(-a+2A),
\end{split}\eee
whose solution is straightforward.
\qed\\

\noindent{\emph{Remark: }}
Note that when we apply the previous result to our case, the functions $I_1[g]$ and $I_2[g]$ will be time depending functions, more precisely
\begin{align*}
I_1[g] &= \int_{-A}^{p_I} f_I\;dx - \int_{0}^{B} f^0\;dx- \int_{0}^{p(t)} \lp f^0 + g \rp\;dx,\\
I_2[g] &= \int_{p_I}^{B} f_I\;dx - \int_{0}^{B} f^0\;dx+ \int_{0}^{p(t)} \lp f^0 + g\rp \;dx,
\end{align*}
where $p_I$ is the unique zero of the initial data $f_I$ and $p(t)$ the unique zero of the function $f(\cdot,t) = f^0(\cdot) + g(\cdot,t)$.
\\
% ====================================================================================================================
% ====================================================================================================================

\section{Non-linear theory}\label{section-non-linear}

\setcounter{equation}{00}

Let $f^0$, $p^0$ be an equilibrium point as specified by the hypothesis in the theorem.
We try to construct a solution of equation
$$f_t=f_{xx}-f_x(p(t))\left[\delta_{p(t)-a}-\delta_{p(t)+a}\right] $$
as a perturbation of this stationary state, i.e.,
$$f(x,t)=f^0(x)+g(x,t).$$

Consider the space $\mathcal V$ of functions $f:[-A,B]\to\mathbb R$ that have only one root $p\in (p^0-\nu,p^0+\nu)$ and that are positive in $[-A,p)$, negative in $(p,B]$. We will allow perturbations  $g$ such that  $g(\cdot,t)$ for each fixed time belongs to the space
\be\label{neighborhood}U_\rho=\left\{g\in Y\,:\;\norm{g}_Y<\rho,\, f^0+g\in \V\right\},\ee
for some $\rho$.
It is essential that this neighborhood is open and uniform in some sense:

\begin{lemma}\label{lemma-neighborhood}
Given $\rho>0$, and $0<\omega\leq\rho$, if $\omega<\min\left\{\lambda^0\nu,\lambda^0\right\}$, the neighborhood
$$\left\{g\in Y\;:\; \norm{g}_Y< \omega\right\}$$
is contained in $U_\rho$.
\end{lemma}

\noindent{\emph{Proof: }}
Assume that $\norm{g}_Y< \omega$ as specified in the hypothesis of the lemma. Let us prove that $f=f^0+g$ has an unique zero in the interval $(-\nu,\nu)$.
First of all note that $f$ is strictly monotone decreasing in the interval $(-\nu,\nu)$ since
$$
f_x (x)= f^0_x(x)+g_x(x)= -\lambda^0 + g_x(x)\le -\lambda^0+\omega<0,$$
for $\omega<\lambda^0$. Moreover
 $$f(-\nu)=|\lambda^0|\nu +g(-\nu)>|\lambda^0|\nu-\omega >0,\quad f(\nu)=-|\lambda^0|\nu +g(\nu)< - |\lambda^0|\nu+\omega<0.$$
 This implies that $f$ has an unique zero in the interval $(-\nu,\nu)$.
 Let us call this unique zero $p$. It holds that $-f_x(p)\geq \lambda^0-w > 0$. A similar argument gives that $f$ is positive in $[-A,p)\cup(p,B]$.
 \qed\\

\medskip

Now we relate the nonlinear problem to the linear one, in order to obtain a formula for the non-linear part. Let $f(x,t)$, $p(t)$ be a solution to our problem. We write $f(x,t)=f^0(x)+ g(x,t)$, and $p(t)=p^0+ q(t)$. Assume that $q\in(-\nu,\nu)$. As in the calculation in Section \ref{section-linearized}, we expand for each time $t$. Thus
\be\label{formula100}0=f(p(t),t)=f^0(p^0)+ q f^0_x(p^0)+g(p^0,t)+R_1(t),\ee
where
$$R_1(t)=g(p(t),t)-g(p^0,t).$$
But $f^0(p^0)=0$ and $f_x^0(p^0)=-\lambda^0$. For simplicity in the notation, we drop the dependence in $t$.
Then from \eqref{formula100} we can write
\be\label{formula-q}
q=\frac{g(p^0)+R_1}{\lambda^0}.\ee
On the other hand, we expand
\be\label{Taylor1}
f_x(p(t),t)=f^0_x(p^0)+ g_x(p^0,t)+ R_2=-\lambda^0+ g_x(p^0,t)+ R_2,\ee
where
$$R_2(t)=g_x(p(t),t)-g_x(p^0,t),$$
and also
\be\label{Taylor2}\delta_{p(t)\pm a}=\delta_{p^0 \pm a}+q\delta'_{p^0\pm a}+R_3^{\pm}.\ee
When we substitute \eqref{formula-q}, \eqref{Taylor1} and \eqref{Taylor2} into the equation for $f$, after some computation we obtain the following equation (where we group the terms of the same  order)
\be\label{expansion-epsilon}\begin{split}
0&=f^0_t-f^0_{xx}+f^0_x(p^0)[\delta_{p^0-a}-\delta_{p^0+a}]\\
&+\lp g_t-Lg\rp \\
&- N(g),
\end{split}\ee
for
\be\label{nonlin-Oper}\begin{split}
N(g) =&\lambda_0 (R_3^- - R_3^+)\\
&+ R_1 \lp \delta'_{p^0-a}-\delta'_{p^0+a}\rp\\
&-g_x(p^0) \left[\lp \delta_{p-a}-\delta_{p^0-a} \rp -\lp \delta_{p+a}-\delta_{p^0+a} \rp \right] \\
& -R_2\lp \delta_{p-a}-\delta_{p+a}\rp \\
=: & N^1(g)+N^2(g)+N^3(g)+N^4(g).
\end{split}\ee
Thus we have proved that

\begin{prop}\label{prop-relation}
Fixed any admissible equilibrium $f^0,p^0$, then problem \eqref{non-linear}-\eqref{initial-data} is equivalent to
\be\left\{\begin{split}
& g_t=Lg+N(g) \\
& g(x,0)=g_I\\
& g_x(-A,0)=g_x(B,0)=0,
\end{split}\right.\ee
for some $g(x,t)$ that satisfies
$$f^0+g(\cdot,t)\in \V,$$
where $L$ and $N$ are given in  \eqref{op_L}, \eqref{nonlin-Oper}, respectively, and $g_I=f_I-f^0$.
The equivalence is simply
$$f=f^0+g.$$
\end{prop}

\bigskip
% =================================================================================================================

\subsection{Lipschitz estimates}

We have that the non-linear part $N$ is locally Lipschitz in $U_\rho$:

\begin{lemma}\label{lemma-Lipschitz0}
Fix $\rho\leq c\;\lambda^0$ for some $c$ positive constant less than one. We have that the operator $N$ given by \eqref{nonlin-Oper} is Lipschitz in the neighborhood $U_\rho$, given in \eqref{neighborhood}. More precisely,
$$\norm{N(g)-N(\tilde g)}_X\leq C\norm{g-\tilde g}_Y,$$
for all $g,\tilde g\in U_\rho$. The constant $C$ depends on $\rho$ and $\lambda_0$.
\end{lemma}

\noindent{\emph{Proof: }}
Let $g,\tilde g\in U_\rho$, with roots $p,\tilde p\in(p^0-\nu,p^0+\nu)$. Let $p=p^0+q$, $\tilde p=p^0+\tilde q$. From \eqref{formula-q} we can estimate the difference $\abs{q-\tilde q}$ by
\bee\begin{split}
\abs{q-\tilde q} &\leq \frac{1}{\lambda^0} \abs{g(p)-\tilde g(\tilde p)} \\
&\leq \frac{1}{\lambda^0} \left\{\abs{g(p)- g(\tilde p)}+\abs{g(\tilde p)-\tilde g(\tilde p)}\right\} \\
& \leq \frac{1}{\lambda^0} \left[ \abs{q-\tilde q}\cdot \sup_{\xi\in(-\nu,\nu)}\left\{\abs{g_x(p^0+\xi)}\right\}+
\sup_{\xi\in(p^0-\nu,p^0+\nu)}\left\{\abs{{g(p^0+\xi)-\tilde g(p^0+\xi)}}\right\}\right].
\end{split}\end{equation*}
If
$$ \frac{\rho}{\lambda^0} = c<1,$$
it holds
$$\abs{q-\tilde q}\leq \frac{1}{ (1-c)\lambda^0}  { \norm{g-\tilde g}}_{\mathcal C^0(p^0-\nu,p^0+\nu)}.$$

Now we compute $\norm{N(g)-N(\tilde g)}_X$. The definition of negative fractional Sobolev norm of order $-\theta$, $\theta>0$, is
$$\norm{\mu}_{H^{-\theta}}=\sup_{\norm \phi_{H^\theta}=1}\langle \phi,\mu\rangle.$$
In our case, we will take $\theta=-r+2\alpha>0$. Note that for an exponent  $3/2 <\theta<2$, we have that $H^{\theta}(-A,B)\subset \mathcal C^{1,\beta}(-A,B)$.\\

Since $R_3^{\pm}=\delta_{p\pm a}-\delta_{p^0\pm a}-q\delta'_{p^0-a}$,
we can easily estimate $\norm{N^1(g)-N^1(\tilde g)}_X$, as follows:
$$\norm{R_3^{\pm} - \tilde R_3^{\pm}}_X\leq \norm{\delta_{p\pm a} -\delta_{\tilde p\pm a}} _{H^{-\theta}(-A,B)}
+\abs{q-\tilde q} \cdot\norm{\delta'_{p^0\pm a}}_{H^{-\theta}(-A,B)}.$$
The first term can be bound by
$$\norm{\delta_{p\pm a} -\delta_{\tilde p\pm a}} _{H^{-\theta}(-A,B)}\leq \sup_{\norm{\phi}_{H^\theta}=1} \left\{\abs{ \phi(p\pm a)-\phi(\tilde p\pm a)}\right\}\leq \abs{p-\tilde p}=\abs{q-\tilde q},$$
taking into account that $H^{\theta-1}\subset L^\infty$ for $3/2 < \theta < 2$.\\

For the estimate of $N^2(g)-N^2(\tilde g)$, it holds
$$\norm{\lp R_1 -\tilde R_1\rp\delta'_{p^0-a}}_X\leq \norm{R_1-\tilde R_1}\norm{\delta'_{p^0-a}}_{H^{-\theta}},$$
where
\bee\begin{split}
\norm{R_1-\tilde R_1} & \leq \abs{g(p)-\tilde g(p)}+\abs{\tilde g(p)-\tilde g(\tilde p)}+\abs{g(p^0)-\tilde g(p^0)} \\
&\leq 2 \norm{g-\tilde g}_{\mathcal C^0(p^0-\nu,p^0+\nu)}+\norm{g}_{\mathcal C^1(p^0-\nu,p^0+\nu)}\abs{q-\tilde q}.
\end{split}\eee
On the other hand, for $N^3$,
\bee\begin{split}
&\norm{g_x(p^0)\left[\delta_{p-a}-\delta_{p^0-a}\right]-\tilde g_x(p^0)\left[\delta_{\tilde p-a}-\delta_{p^0-a}\right]}_X \\
& \leq \abs{g(p^0)-\tilde g(p^0)}\norm{\delta_{p-a}-\delta_{p^0-a}}_{H^{-\theta}}+\abs{\tilde g_x(p^0)}\norm{\delta_{p-a}-\delta_{\tilde p-a}}_{H^{-\theta}} \\
&\leq \norm{g-\tilde g}_{\mathcal C^0(p^0-\nu,p^0+\nu)} +\norm{\tilde g}_{\mathcal C^1(p^0-\nu,p^0+\nu)}\norm{\delta_{p-a}-\delta_{\tilde p-a}}_{H^{-\theta}}.
\end{split}\eee
And for the last term $N^4$,
\bee\begin{split}\norm{ R_2\delta_{p-a}-\tilde R_2\delta_{\tilde p-a}}_X &
\leq\norm{R_2\lp\delta_{p-a}-\delta_{\tilde p-a}\rp}_X+\norm{\lp R_2-\tilde R_2\rp \delta_{\tilde p-a}}_X \\
& \leq \abs{R_2}\norm{\delta_{p-a}-\delta_{\tilde p-a}}_{H^{-\theta}}+\abs{ R_2-\tilde R_2}\norm{\delta_{\tilde p -a}}_{H^{-\theta}}.
\end{split}\eee
Since $$\abs{R_2}=\abs{g_x(p)-g_x(p^0)}\leq \norm{g}_{\mathcal C^2(p^0-\nu,p^0+\nu)}\abs{q},$$
it holds
\be\label{estR2}\begin{split}
\abs{R_2-\tilde R_2} & \leq \abs{g_x(p)-g_x(\tilde p)}+\abs{g_x(\tilde p)-\tilde g_x(\tilde p)}+\abs{\tilde g_x(p^0)-g_x(p^0)} \\
& \leq \norm{g}_{\mathcal C^2(p^0-\nu,p^0+\nu)}\abs{q-\tilde q}+2\norm{g-\tilde g}_{\mathcal C^1(p^0-\nu,p^0+\nu)}\\
& \le \lp \frac{\rho}{\lambda^0(1-c)}+2\rp \norm{g-\tilde g}_Y.
\end{split}\ee
In particular, there exists a constant $C$,  such that
$$\norm{N(g)-N(\tilde g)}_X\leq C \norm{g-\tilde g}_Y,$$
for every $g$, $\tilde g\in U_\rho$, where $C = \textrm{max}\left\{1,\frac{\rho}{\lambda^0(1-c)}\right\}$ .
\qed\\

\bigskip

The previous estimate can be improved near the origin, in fact, $DN(0)=0$. This estimate is needed in order to pass from local to global existence. More precisely,

\begin{lemma}\label{lemma-growth}
We have that
\be\label{growth}\norm {N(g)}_X= o \lp\norm{g}_{Y}\rp,\mbox{ when } \norm{g}_{Y}\to 0.\ee
\end{lemma}

\noindent{\emph{Proof: }}
The previous lemma implies that $\norm {N(g)}_X= O \lp\norm{g}_{Y}\rp$ when $\norm{g}_{Y}\to 0$. For the improvement we need, take a perturbation of the form $f=f^0+\epsilon g$ for some $\epsilon>0$ small. Let $p$ be the root of $f$ near $p_0$. Because of \eqref{formula-q} we can write $p=p^0+\epsilon q$ for $q=\frac{g(p)}{\lambda_0}$. Assume $\norm{g}_Y=1$, for simplicity.
First, note that because of the Taylor theorem we can write
$$R_3^{\pm}=\epsilon q\left[\delta'_{x^{\pm}_3\pm a}-\delta'_{p^0\pm a}\right],$$
for some $\abs{x_3-p^0}\leq \epsilon \nu$. We compute
\be\label{delta'}\norm{\delta'_{x_3-a}-\delta'_{p^0-a}}_{H^{-\theta}}=\sup_{\norm \phi_{H^\theta}=1} \abs{\phi'(x_3-a)-\phi'(p^0-a)},\ee
where $\theta = -r+2\alpha$. Then, because of our choice of $\theta$, $H^\theta(-A,B)\subset \mathcal C^{1,\beta}(-A,B)$, and from \eqref{delta'} we deduce that
$$\norm{\delta'_{x_3-a}-\delta'_{p^0-a}}_{H^{-\theta}}\leq c\epsilon^\beta.$$
As a consequence,
$$\norm{N^1(\epsilon g)}_X\leq C \epsilon^{1+\beta}.$$
Next, note that $R_1=\epsilon \left[ g(p)-g(p^0)\right]$, so that $\abs{R_1}\leq \epsilon^2 \norm{g}_{\mathcal C^1} q$, and this gives a bound for $N^2$.
We handle the rest of the terms in a similar manner, and we conclude that
$$\norm {N(g)}_X\leq C_2 \epsilon^{1+\beta}.$$
for some $C_2$ depending on $\frac{1}{(1-c)\lambda^0}$, but independent of $\norm{g}_Y$.
The lemma is proved.
\qed\\

\noindent{\emph{Remark: }}
The choice of $s$ in the definition of the function space $Y$ is made here. In particular, we need $\mathcal C^{2,\beta}$ regularity near the origin in order to have a Lipschitz estimate in Lemma \ref{lemma-Lipschitz0}.
\\

% ========================================================================================================
% ========================================================================================================

\subsection{Local existence}

First we review some basic facts about the nonlinear evolution equation
\be\label{non-linear-problem}\left\{\begin{split}
&g_t=Lg+N(g), \\
&g(0)=g_I.
\end{split}\right.\ee
Here we always assume that $L$ is a sectorial operator, $U$ an open subset of $Y=X^\alpha$ for some $0\leq \alpha<1$ and $N:U\to X$ a  locally Lipschitz function. More precisely, given any $g\in U$, there exists a neighborhood $V\subset U$ of $g$ such that for all $g,\tilde g\in V$, it holds
$$\norm{N(g)-N(\tilde g)}\leq C\norm{g-\tilde g}_Y,$$
for some constant $C$ depending on the neighborhood $V$.\\

A \emph{solution} of this initial value problem on $(0,t_1)$ is a continuous function $g:[0,t_1)\to X$ such that $g(0)=g_I$; on $(0,t_1)$ we have $g(t)\in U$, $g(t)\in D(L)=Z$, $\frac{dg}{dt}$ exists, $t\to N(g(t))$ is locally H\"older continuous, $\int_{0}^{0+\mu}\norm{N(g(t))}dt<\infty$ for some $\mu>0$; and moreover the differential equation \eqref{non-linear-problem} is satisfied on $(0,t_1)$. In particular, we have the Duhamel's formula
$$g(t)=e^{Lt}g_I+\int_{0}^t e^{L(t-s)}N(g(s))ds.$$

It is well known that

\begin{thm}[Theorem 3.3.3 in \cite{Henry:parabolic}]\label{thm1}
Assume that $L$ is a sectorial operator, $0\leq \alpha<1$, and $N:U\to X$, $U$ an open subset of $X^\alpha$, and $N$ is locally Lipschitz. Then for any $g_I\in U$, there exists $0<T=T(g_I)$ such that \eqref{non-linear-problem} has a unique solution $g$ on $(0,T)$ with initial value $g(0)=g_I$.
\end{thm}

\noindent{\emph{Remark: }}
In fact, it is true that
$$g\in\mathcal C^1 \lp[0,T),X\rp.$$
\\

Let $g_I=f_I-g^0$. The previous theorem assures local existence in time for $g$ if $g_I\in U_\rho$. Then, Proposition \ref{prop-relation} gives that,
the function $f(x,t)=f^0(x)+g(x,t)$ we construct is a solution to our original problem. Moreover, it has the right signs: it has a unique root in $[-A,B]$, call it $p(t)$, and it is positive in $[-A,p(t))$, negative in $(p(t),B]$. The proof is similar as in Lemma \ref{lemma-neighborhood}.

\bigskip

% ====================================================================================================================
% ====================================================================================================================

\subsection{Existence of a center manifold}

In the present section we give the proof the existence of a center manifold for our problem. Here we use follow the approach of chapter 6 in the book \cite{Henry:parabolic}, that we find very clear and suitable for our purposes;  although \cite{Chen-Hale-Tan:invariant-foliations} is the classical one for applications in partial differential equations. Other references are, for instance, \cite{Chow-Hale:bifurcation-theory} and \cite{Carr:centre-manifold}.

We consider the equation \eqref{non-linear-problem}. A set $S\subset X^\alpha$ is a \emph{local invariant manifold} if for any $g_I\in S$, there exists a solution $g(\cdot)$ of the differential equation on an open interval $(t_1,t_2)$ containing $0$ with $g(0)=g_I$ and $g(\cdot)\in S$ for $t_1<t<t_2$. We say that $S$ is an \emph{invariant manifold} if we can always choose $(t_1,t_2)=(-\infty,+\infty)$.

In particular, for a coupled system
\bee
\left\{
\begin{split}
&(g_1)_t=N_1(g_1,g_2), \\
&(g_2)_t=L_2 (g_1+g_2)+N_2(g_1,g_2),
\end{split}\right.
\eee
we look for invariant manifolds of the form $S=\left\{(g_1,g_2)\,:\, g_2=\sigma(g_1)\right\}$. It is well known that:

\begin{thm}[Theorem 6.1.2. in \cite{Henry:parabolic}]\label{thm-existence-manifold}
Let $X_1,X_2$ be Banach spaces and assume $-L_2$ is sectorial in $X_2$. Let $U$ be a neighborhood of the origin in $X_2^\alpha$ for some $\alpha<1$. Suppose that
\begin{enumerate}
\item The operator $N_2:X_1\times U\to X_2$ is locally Lipschitz with
\be\label{lipschitz2}\norm{N_2(g_1,g_2)-N_2(\tilde g_1, \tilde g_2)}_{X_2}\leq \kappa \lp \norm{g_1-\tilde g_1}_{X_1}+\norm{g_2-\tilde g_2}_{X_2^\alpha}\rp,\ee
and
\be\label{formula-N}\norm{N_2(g_1,g_2)}_{X}\leq N.\ee
\item The operator $N_1: X_1\times U\to X_1$ is locally Lipschitz with
\be\label{lipschitz1}\norm{N_1(g_1,g_2)-N_1(\tilde g_1, \tilde g_2)}_{X_1}\leq \mu\norm{g_1-\tilde g_1}_{X_1}+ M_2\norm{g_2-\tilde g_2}_{X_2^\alpha}.\ee
\item There exists some $\gamma_0>0$ such that
$$\norm{e^{L_2 t}}_{X_2}\leq C_1 e^{-\gamma_0 t},\quad \norm{(-L_2)^{\alpha}e^{L_2 t}}_{X_2}\leq C_1 t^{-\alpha}e^{-\gamma_0 t},\quad \mbox{for }t>0.$$
\item For some positive constant $D$,
$$\left\{g_2 : \norm{g_2}_{X_2^\alpha} \leq D\right\}\subset U\quad \mbox{and}\quad C_1N\int_{0}^{\infty} t^{-\alpha} e^{-\gamma_0 t}dt<D.$$
\item For some positive constant $E$,
\be\label{theta}\theta:=\kappa C_1 \int_{0}^\infty u^{-\alpha}e^{-\gamma_0 u} e^{(\mu+EM_2)u}du\ee
has $$\theta\leq \frac{E}{1+E} \texto{and }\theta\max\left\{ 1,\frac{(1+E)M_2}{\mu+M_2 E}\right\}<1.$$
\end{enumerate}
Then there exists an invariant manifold (the center manifold)
$$S=\left\{ (g_1,g_2) : g_2=\sigma(g_1), -\infty<t<+\infty, g_1\in X_1\right\},$$
with
$$\norm{\sigma(g_1)}_{X_2^\alpha}\leq D,$$
and
$$\norm{\sigma(g_1)-\sigma(\tilde g_1)}_{X_2^\alpha}\leq E \norm{g_1-\tilde g_1}_{X_1}.$$
\end{thm}

\noindent{\emph{Remark: }}
If  $\mu<\gamma_0$ and $\kappa$ is sufficiently small, we can easily satisfy condition \emph{5} above.
\\
\noindent{\emph{Remark: }}
Note that condition \emph{4} and \emph{5} ensure global in time existence of solutions, that have been proven to exist only locally in time in Theorem \ref{thm1}. In fact, if the constants $C_1$, $C_2$,, $N$ etc satisfy the above conditions, it is easy to check that (see Theorem 5.1.1 in \cite{Henry:parabolic})
$$
\norm{g_2}_{X_2^\alpha} <\rho,\quad \forall t\in (0, T).
$$
This implies that the solution exists for all time (Theorem 3.3.4 in \cite{Henry:parabolic}) and stays in the neighborhood $U$, i.e.,
$$\norm{g_2}_{X_2^\alpha} <\rho,\quad \forall t\in (0, +\infty).$$\\

In the following, we check that Proposition \ref{thm-existence-manifold} can be applied in our problem.  First, Proposition \ref{prop-basic} shows that we can split $X=X_1\oplus X_2$ with $X_1=\ker (L)$. The subspace $X_1$ has finite dimension two. We write $g=g_1+g_2$, $g_1\in X_1$ and $g_2\in X_2$.
Let $L_i$ be the restriction of the linear operator $L$ to $X_i$, $i=1,2$, and let $P$ be the projection on $X_1$ along $X_2$. We denote
$$N_1(g_1,g_2):=PN(g_1+g_2),\quad N_2(g_1,g_2):= (I-P)N(g_1+g_2).$$
Thus, our problem \eqref{non-linear-problem} is equivalent to
\be\label{system}
\left\{
\begin{split}
&(g_1)_t=N_1(g_1,g_2), \\
& (g_2)_t=L_2(g_2)+N_2(g_1,g_2).
\end{split}\right.
\ee
In order to apply Theorem \ref{thm-existence-manifold}, we need to check the Lipschitz condition of the non-linear terms.  In the following proposition we give a characterization of $N_1$ and $N_2$.

\begin{prop}\label{lemma-formula-N1}
Consider any $g(\cdot,t)\in X$ such that $f^0+g(\cdot, t)$ has a unique root in $[-A,B]$, denote it by  $p(t)$. Then there exist constants $\gamma_1$ and $\gamma_2$ depending only on $A,B,a$ such that if $g$ is written as $g=g_1+g_2$, $g_i\in X_i$, $i=1,2$, we can write
$$N_1(g_1,g_2)=\gamma_1 R_2(t) g_0+\gamma_2 R_2(t) h_0,$$
where $g_0$, $h_0$ are the eigenfunctions of the zero eigenvalue given in Proposition \ref{prop-eigenfunctions}, and $R_2(t) $ is defined as in \eqref{Taylor1}
$$R_2(t):=g_x(p(t))-g_x(p_0).$$
\end{prop}

\noindent{\emph{Proof: }}
We write
\be\label{N_1'}N(g)=N_1(g_1,g_2)+N_2(g_1,g_2),\ee
where $N_i(g_1,g_2)\in X_i$, $i=1,2$. Because $\{g_0,h_0\}$ is a basis for $X_1$, we can write
\be\label{N_1}N_1(g_1,g_2)=k(t) g_0+ z(t) h_0,\quad \int_{-A}^0 N_2( g_1,g_2) \;dx  = \int_{0}^B N_2( g_1,g_2) \;dx =0, \ee
for some $k(t)$, $z(t)$ functions. The exact value of these coefficients can be computed thanks to Lemma \ref{characterize-kernel}.

On the other hand, consider the expression for $N(g)$ given in \eqref{nonlin-Oper}. If we integrate it from $0$ to $B$, for each fixed time, we obtain that most of the terms cancel and we are left with
$$I_2[N(g)]:=\int_0^B N(g)\,dx  =R_2(t).$$
A similar integration from $-A$ to $0$ gives that
$$I_1[N(g)]:=\int_{-A}^0 N(g)\,dx=- R_2(t),$$
where $R_2(t)$ is defined as in \eqref{Taylor1}. It follows that $I_1[N(g)] = -I_2[N(g)]$ and from \eqref{coefficients} we get
$$k(t)=\frac{2(-a+A+B)}{a(a-2A)(a-2B)}R_2(t),\quad  z(t)=\frac{A-B}{(a-2A)(a-2B)}R_2(t),$$
as desired.
\qed\\

With all these ingredients we can prove that:

\begin{prop} \label{prop-center-manifold}
Fixed an equilibrium point $f^0$, there exists some $\rho >0$ such that there exists a center manifold for \eqref{system} in the neighborhood $U_\rho$ with
$$\norm{\sigma(g_1)}_{\alpha} \leq D,$$
and Lipschitz condition
$$\norm{\sigma(g_1)-\sigma(\tilde g_1)}\leq E \norm{g_1-\tilde g_1}.$$
for $D$, $E$ given as in Theorem \ref{thm-existence-manifold}.
\end{prop}

\noindent{\emph{Proof: }}
It is a consequence of Theorem \ref{thm-existence-manifold}, but we need to check carefully the constants involved. First, note that the operator $-L_2$ generates an analytic semigroup and satisfies the estimates given in \eqref{estimates-semigroup}
$$\norm{e^{L_2 t}}_{X_2}\leq C_1 e^{-\gamma_0 t},\quad \norm{(-L_2)^{\alpha}e^{L_2 t}}_{X_2}\leq C_1 t^{-\alpha}e^{-\gamma_0 t},$$
for some $0<\alpha<1$ and $0<\gamma_0<\hat \gamma$.

Next, we look at the non-linear part. We will work in a neighborhood $U_\rho$ for some $\rho<\min\{\lambda_0,\nu\lambda_0\}$ so that we are in the situation of Lemma \ref{lemma-neighborhood}.

Lipschitz condition \eqref{lipschitz2} follows from Lemma \ref{lemma-Lipschitz0}. On the other hand, \eqref{lipschitz1} is a direct consequence of Proposition \ref{lemma-formula-N1} and (\ref{estR2}). The Lipschitz constants $\kappa,\mu,M_2$ depend on $\rho$. More precisely,
these constants can be estimated by
$$\frac{1}{(1-c)\lambda^0} O(\rho)\quad \mbox{in }U_\rho,$$
when $\rho\leq c\lambda^0$ for some $c<1$. Fix $\rho_5$ such that condition {\emph 5} is satisfied for all $\rho\leq \rho_5$.

Now, in order to check condition \emph{4}, we need Lemma \ref{lemma-growth}, since the Lipschitz estimate for $N$ is not enough. We have seen that in this case,
$$\norm{N(g)}_X = C_2\norm{g}_Y^{1+\beta}, \quad \mbox{when}\norm{g}_Y\to 0,$$
for some $0<\beta<1$ and $C_2=O\lp\frac{1}{(1-c)\lambda^0}\rp$. Then looking at \eqref{formula-N} we can write
$$N= \frac{2}{\lambda^0}  \rho^{1+\beta} \quad\mbox{in }U_\rho.$$
Choose $D= c \rho$ , with $c<1$, we get
$$C_1\frac{2}{\lambda^0}  \rho^\beta\int_0^\infty u^{-\alpha} e^{-\gamma_0 u}du<c.$$
Choosing $\rho$ small enough, condition 4. is satisfied.

We conclude that  the solution exists globally in time and there exists an invariant manifold
$$S=\left\{ (g_1,g_2) : g_2=\sigma(g_1), -\infty<t<+\infty, g_1\in X_1\right\},$$
as we wished.
\qed\\

\bigskip

% ==============================================================================================
% ==============================================================================================

\subsection{Stability of a family of equilibria}

Now we try to understand the decay near the center manifold. We are given an initial condition $f_I$ with masses
$$\int_{-A}^{p_I} f_I=m_1, \quad -\int_{p_I}^B f_I =m_2.$$
Assume that it satisfies $\norm{f_I-f^0}_Y<\rho$ for some $\rho$ small enough that will be specified later.
Here $f^0$ is any of the stationary states described in Section \ref{section-stationary} with masses $\tilde m_1,\tilde m_2$ (not necessarily equal to $m_1,m_2$). We can take, without loss of generality, $p^0=0$.

By linearizing around $f^0$, we are going to show that there exists a solution $f$ of the system with initial condition $f_I$, that exists for all time, and that decays exponentially to a stationary state $f_\infty$. Moreover,  $f_\infty$ is chosen as the unique stationary state with
$$\int_{-A}^{p_\infty} f_\infty=m_1, \quad -\int_{p_\infty}^B f_\infty =m_2, $$
for $f_\infty(p_\infty)=0$.\\

Global existence and asymptotic decay are given by

\begin{thm}[Theorem 6.1.4 and Corollary 6.1.5. in \cite{Henry:parabolic}]\label{thm-Henry-decay}
In addition to the hypothesis of Theorem \ref{thm-existence-manifold}, assume that  $\dim X_1<\infty$, and
\be\label{formula-r}
r:=\theta\lp 1+M_2\frac{1+E}{\gamma_0-\mu'}\rp<1,\ee
for $\mu'=\mu+M_2 E$.
Then $S$ is uniformly asymptotically stable with asymptotic phase. Specifically, there exist $\delta>0$, $c>0$ so that any solution $(g_1(t),g_2(t))$ with initial condition satisfying
$$\norm{g_2(0)-\sigma(g_1(0))}_{X_2^\alpha}<\delta,$$
exists for all time $t\geq 0$ and there is a solution $\bar g_1(t)$ of
$$(\bar g_1)_t=N_1(\bar g_1,\sigma(\bar g_1))$$
such that, for $t\geq 0$,
$$\norm{g_1(t)-\bar g_1(t)}+\norm{g_2(t)-\sigma(\bar g_1(t))}_{X_2^\alpha} \leq
C e^{-\gamma(t)}\norm{g_2(0)-\sigma(g_1(0))}_{X_2^\alpha},$$
where
$$\gamma=\gamma_0-(\gamma-\mu')r^{\frac{1}{1-\alpha}}>0.$$
\end{thm}

\noindent{\emph{Remark: }}
In the proof of this theorem, $\delta$ is chosen so that the solution stays in the domain of existence for all time $t\geq 0$, that it is equivalent to
\be\label{delta}kC_1e^{-\gamma(t)}\delta <D/2.\ee
Moreover,
$$c=\frac{M_2(C_1K)^2}{\gamma-\mu'},$$
where $k$ and $K$ are a constants depending only on $\alpha$.

\medskip

Thus in our problem,

\begin{cor}
The center manifold constructed in Proposition \ref{prop-center-manifold} is an attractor if we start with initial condition $g_I\in U_\rho$ for some $\rho$ small enough.
\end{cor}

\noindent{\emph{Proof: }}
We take $\rho$ as in the proof of Proposition \ref{prop-center-manifold}. The value of $\theta$ is given in \eqref{theta}. It is clear that we can take a smaller $\rho$ so that \eqref{formula-r} is satisfied and $\rho\leq \delta$ for  $\delta$ given in \eqref{delta}. The corollary is a consequence of Theorem \ref{thm-Henry-decay} and the remark afterwards.
\qed\\

\bigskip

Up to this point, we have reduced the problem to understand how the trajectories look like inside the center manifold.

\begin{prop}
Any trajectory in the center manifold of Proposition \ref{prop-center-manifold} must be stationary.
\end{prop}

\noindent{\emph{Proof: }}
A trajectory in the center manifold is given by
$$g=g_1+g_2, \quad g_1=c(t)g_0+d(t)h_0,\quad g_2=\sigma(g_1).$$
In order to find the coefficients $c(t),d(t)$, we substitute this solution in \eqref{system},
$$\dot{c} g_0+\dot{d} h_0+\partial_t\lp\sigma(g_1)\rp=L(g_1+g_2)+N(g_1+g_2),$$
and we project onto $X_1$, thus
$$\dot{c} g_0+\dot{d}h_0=N_1(g_1+g_2).$$
We have an explicit formula for the right hand side, thanks to Proposition \ref{lemma-formula-N1}. Then, looking at each eigenspace, we obtain a system of ordinary differential equations for $c(t),d(t)$:
\be\label{reduced-ODE}\left\{\begin{split}
\dot{c}=\gamma_1 R_2(t), \\
\dot{d}=\gamma_2 R_2(t),
\end{split}\right.\ee
where $\gamma_1,\gamma_2$ are explicit constants depending only on $a,A,B$.

We claim that this system of ODE's has a unique global solution. For this, consider
$$
g(x,t) = c(t) g_0 + d(t) h_0+\sigma(c(t) g_0 + d(t) h_0).
$$
 It is easy to check that, for $-\nu\le p(t)\le \nu$ and $g_0$, $h_0$ as in Proposition \ref{prop-eigenfunctions}, it holds
$$
g_x(p(t),t) = c(t)( 1+ \sigma'_{|_{x=p(t)}}),\quad \quad g_x(0,t) = c(t)( 1+ \sigma'_{|_{x=0}}).
$$
Taking into account that $R_2(t) = g_x(p(t),t) - g_x(0,t)$, the above system becomes
\be\label{reduced-ODEII}\left\{\begin{split}
\dot{c}=\gamma_1 c(t)(  \sigma'_{|_{x=p(t)}}-  \sigma'_{|_{x=0}}), \\
\dot{d}=\gamma_2 d(t)(  \sigma'_{|_{x=p(t)}}-  \sigma'_{|_{x=0}}).
\end{split}\right.\ee
Theorem \ref{thm-existence-manifold} ensures that $\sigma$ is a Lipschitz function. It follows that the right-hand side of (\ref{reduced-ODEII}) is also a Lipschitz function with respect to the function $c(t)$.  Existence of a unique solution locally in time is proven.

We claim now that $c(t) = c_I$ and $d(t) = d_I$ solve system \eqref{reduced-ODEII}. Let $g_1(x,t) = c_I g_0(x) + d_I h_0(x)$. We integrate for $x\in(-A,0)$  equation
$$
(g_1)_t = N_1(g_1,g_2),
$$
we obtain that
$$
0 = \frac{d}{dt} \int_{-A}^0 g_1\;dx = \int _{-A}^0 N_1\;dx =R_2(t).
$$
Thus, the claim is proved. This implies that system \eqref{reduced-ODE} has an unique global in time solution $c(t) = c_I$, $d(t) = d_I$ (the constant one).
\qed\\

\bigskip

In general, the center manifold is not unique and it does not consist only of the stationary points. However, in our case, the set of stationary states has dimension two and is contained in the center manifold. This is a special case and appears in other works, see for instance \cite{Gallay-Wayne:dim-2} on the Navier-Stokes equation, that has a similar behavior. In particular,

\begin{cor}
Our center manifold is exactly the set of equilibrium points for \eqref{system}, i.e.
$$\left\{g\in U_\rho \;:\; g=g_1+g_2, g_1\in X_1, g_2=\sigma(g_1) \right\}=\left\{g\in U_\rho\;:\;Lg + N(g)=0 \right\}.$$
\end{cor}

Now we can finish the proof of Theorem \ref{theorem-center-manifold}. For each initial datum $f_I$ in the neighborhood $U_\rho$, we have built a solution of the nonlinear problem for all time $t>0$, of the form $f=f^0+g$, and that converges exponentially to an stationary state, that may not be the same as $f^0$. Since we have already classified all the stationary solutions, in particular, because of conservation of mass, it has to be the precise $f_\infty$.

Note that any stationary solution $f_\infty$ can be written as a {\em{shifted}} zero eigenvector in the following way:
\bee f_\infty(x) = \left\{\begin{split}
\lambda_{\infty} g_0(x-p_\infty)\quad x\in[-A+p_\infty, B], \\
\lambda_{\infty} g_0(-A)\quad x\in[-A, -A+p_\infty].
\end{split}\right.\eee

% ===================================================================================================================

\subsection{Some final remarks}

In this last section, we group together a few remarks. First, we can write a more explicit formula for the center manifold. Consider $f^0$ a general equilibrium point with $f^0(0)=0$. Consider $f_\infty$ the unique stationary state corresponding to the initial datum $f^0+g_I$, with $g_I\in U\rho$. Let $g_I$ be decomposed into our basis of eigenfunctions
$$
g_I= c_Ig_0 + d_I h_0 +\tilde{g}_I,\quad   \tilde{g}_I\in X_2.
$$
The corresponding solution $f(x,t) = f^0 + g_1(x,t) + g_2(x,t)$ converges exponentially fast to the equilibrium point
$$
f_\infty = f^0 +   c_Ig_0 + d_I h_0  + \sigma(c_Ig_0 + d_I h_0),
$$
where $c_Ig_0 + d_I h_0+ \sigma(c_Ig_0 + d_I h_0)$ belongs to the center manifold. It holds that
\bee\begin{split}
h_1(c_I,d_I):=&\int_{-A}^0 f_\infty \,dx=\lambda_\infty a \lp p_\infty-\tfrac{a}{2}+A\rp-\tfrac{\lambda_\infty}{2}p_\infty^2= (c_I+\lambda^0)I_1[g_0] + d_I I_1[h_0],\\
h_2(c_I,d_I):=&\int_{0}^B f_\infty \,dx= \tfrac{\lambda_\infty}{2}p_\infty^2-\lambda_\infty a\lp B-p_\infty-\tfrac{a}{2}\rp = (c_I+\lambda^0)I_2[g_0] + d_I I_2[h_0],
\end{split}
\eee
taking into account that any general equilibrium point $f^0$ can be written as $f^0 = \lambda^0 g_0$. Inverse function theorem shows that given $c_I$ and $d_I$, one can uniquely determine $f_\infty$ (or, equivalently, determine $p_\infty$ and $\lambda_\infty$) from the above formula. More precisely
$$
(\lambda_\infty, p_\infty) = H^{-1}(h_1(c_I,d_I), h_2(c_I,d_I)),
$$
where
\bee
H(\lambda,p)=\left(
\begin{array}{c}
\lambda a \lp p-\tfrac{a}{2}+A\rp-\tfrac{\lambda}{2}p^2 \\
\tfrac{\lambda}{2}p^2-\lambda a\lp B-p-\tfrac{a}{2}\rp
\end{array}
\right).\eee
In fact, we can compute its derivative $DH(\lambda^0,0)$ and get
\bee
DH(\lambda^0,0)=\left(
\begin{array}{cc}
a\lp -\tfrac{a}{2}+A\rp & \lambda^0 a\\
-a\lp -\tfrac{a}{2}+B\rp & \lambda^0 a
\end{array}
\right),\eee
that is positive definite because
$$\det DH(\lambda^0,0)=\lambda^0 a^2(B+A-a)>0.$$
The inverse function theorem states the existence of a neighborhood near $(\lambda^0,0)$ such that the function $H$ is invertible (and continuous).
Therefore one can write our function $\sigma$ as
$$
\sigma(c_Ig_0 +d_I h_0) = -(c_I+\lambda^0) g_0 + d_I h_0 + f_\infty ( H^{-1}(h_1(c_I,d_I), h_2(c_I,d_I))).
$$

\bigskip

\noindent{\emph{Remark: }}
The regularity of the free boundary $p(t)$ is given by the following argument: since $f(\cdot,t) \in C^{2,\beta} (-\nu,\nu)$ and $f_x(x,t)\ge c \lambda^0 $ for all $x\in (-\nu,\nu)$, the function
$$p'(t)=-\frac{f_{xx}(p(t))}{f_x(p(t))}$$
is uniformly bounded in time. This implies $p(t)$ is Lipschitz regular $\forall t>0$.\\

\bigskip

And to finish, we would like to see that the size of the neighborhood does not depend on the equilibrium we start with. This is achieved through a careful control of constants. In fact taking as family of equilibria all the functions $f^0$ such that
$$A:= \left\{f^0 \mbox{admissible}\;:\; \lambda^0\geq \chi\right\},$$
we can see that the constants $C$ and $C_2$ in Lemma \ref{lemma-Lipschitz0} and \ref{lemma-growth} are uniformly bounded with respect to $\lambda^0$. Tracing back the constants in the proofs of this section, we have Theorem \ref{thm-independent}.

\bigskip

% ==============================================================================================
% ==============================================================================================

\section{Appendix}
In this last section we give the proof of Lemma \ref{lemma-bounded-operator}. The proof of this result in a more general setting is included in the forthcoming paper \cite{interpolation}. For completeness, we briefly summarize here the main steps. \\
Assume, without loss of generality, that $A=B=1$, $a=3/4$.

Consider the Hilbert space $H^{p}(\T)$ for any $p\in\mathbb R$, with norm given by $$\norm{f}_{H^p}:=\sum (1+n^2)^p|\hat f_n|^2,$$
where $\hat f_n$ are the Fourier coefficients of $f$ in $[-1,1]$, i.e., $f\sim \sum \hat f_n e_n$.

The second functional space we deal with is $H^p_\phi(\T)$, which is defined through the norm
$$\norm{f}_{H_\phi^p}:=\norm{\phi f}_{H^p}.$$
Here the function $\phi$ denotes a smooth cutoff function with support on $[-1/2,1/2]$, identically one on $[-1/4,1/4]$, strictly increasing on $[-1/2,-1/4]$, and strictly decreasing on $[1/4,1/2]$.
The intersection space $X:=H^p\cap H^q_\phi$ has a norm defined as
$$\norm{f}_X:=\norm{f}_{H^p}+\norm{\phi f}_{H^q}.$$

The Laplacian operator acting on $X$ can be written as,
$$\Delta f=\sum \lambda_n \hat f_n e_n, \quad \lambda_n\sim -n^2,$$
where $\lambda_n$ are the eigenvalues of the Laplacian.

Our aim here is to find a functional space $Z$ dense in $X$ such that $\Delta:Z\to X$ is a bounded operator.
We claim that this functional space is given by
 $$Z=H^{p+2}\cap H^{q+2}_{\hat \phi},$$ where
$\hat \phi$ is another cutoff function such that $\hat\phi\equiv 1$ on the support of $\phi$.

The choice of the cutoff functions $\phi$ and $\hat \phi$ implies that
$$\phi \Delta(\hat \phi f)=\phi \Delta(f).$$
As a consequence, the operator $\Delta:Z\to X $ is bounded.

The proof of the density of $Z$ in $X$ requires some technicalities: first, we will prove that periodic $\mathcal C^\infty$-functions are dense in $X$.
 This means that given $f_0\in X$ and $\epsilon>0$, we will find a smooth sequence $f_\epsilon$ such that $\norm{f_\epsilon-f_0}_X<\epsilon$. \\
By taking a collection $\mathcal C$ of overlapping intervals covering $[-1,1]$, we construct a good approximation for $f_0$ in each subintervals of $\mathcal C$. Consider for example $$\mathcal C= \{[-1,-1/2], [-5/8, -1/4], [-3/8, 3/8],[1/4,5/8],[-1/2,1] \} .$$
Once  convergence in each single subinterval is shown, convergence in the whole interval $[-1,1]$ follows from the convergence in the single subintervals by a partition of unity argument.\\
The convergence in $H^p$ is guaranteed by the regularity of the heat equation. Indeed, let $f(x,t)$ be the solution of the heat equation $f_t=\Delta f$ with initial condition $f_0$. Then, if $f_0=\sum \hat f_n e_n$, it is well known that
$$f(x,t)=\sum e^{\lambda_n t} \hat f_n e_n$$
satisfies $f \to f_0$ in $H^p$ when $t\to 0$. Moreover, also $f\to f_0$ in $H^q$ on $[-3/8,3/8]$, away from the vanishing points of the cutoff $\phi$.\\
On the subinterval $[-5/8, -1/4]$ we regularize the function $f_0$ as $f_\epsilon=\theta_\epsilon * f_0$ where $\theta$ is a mollifier with support in $(0,1)$ and $\theta_\epsilon(x)=\frac{1}{\epsilon}\theta(x/\epsilon)$, with supp $\theta_\epsilon =(0,\epsilon)$. Since $f_0\in H^p$, it is known that $f_\epsilon\to f_0$ in $H^p$.

It remains to prove that $\phi f_\epsilon$ converges to $\phi f_0$ in $ H^q$: since $\phi f_0\in H^q$, one can consider a small shifting $\phi(x-t)f_0(x)\in H^q$ for $t\in (0,\epsilon)$, so that $\phi(x)f_0(x+t)\in H^q$. In this way, the function $\phi(x)f_\epsilon(x)$, defined as
\be\label{estimate1}
\phi(x)f_\epsilon(x) =\phi(x)\int_0^\epsilon f_0(x+t)\theta_\epsilon(t)dt,
\ee
still belongs to $H^q$ with norm uniformly bounded with respect to $\epsilon$, depending on $\norm{\phi f_0}_{H^q}$.\\
The term $\norm{\phi f_\epsilon-\phi f_0}_{H^q}$ can be rewritten as
$$\norm{\phi f_\epsilon-\phi f_0}_{H^q}\leq \norm{\phi (f_\epsilon-\psi_\epsilon)}_{H^q}+\norm{\phi (\psi_\epsilon- \psi)}_{H^q}+\norm{\phi (f_0 -\psi)}_{H^q},$$
where  $\psi_\epsilon:=\theta_\epsilon *\psi$. \\
The first term of the sum can be uniformly bounded by $\norm{\phi (f_0-\psi)}_{H^q}<\delta/3$ (in the same way \eqref{estimate1} is bounded by $\norm{\phi f_0}_{H^q}$).\\
Due to a standard property of mollifiers, we have that $\norm{\phi (\psi_\epsilon- \psi)}_{H^q}\le \delta/3$.\\
Finally, for all $\delta >0$, there exists a function $\psi\in \mathcal C^\infty$ with compact support in the interval $[-5/8, -1/4]$ such that $\norm{\phi (f_0 -\psi)}_{H^q}\leq \delta/3$.\\
This implies that for every $\epsilon$ small enough, one can find a $\delta$ such that $\norm{\phi f_\epsilon-\phi f_0}_{H^q}\leq \delta$.

Once the proposition is true for the $\Delta$ operator, then it is true in our case, where we deal with a bounded perturbation of the Laplacian:
$$Lg:=g_{xx}- g_x(0)\left[\delta_{x=-\frac{3}{4}}-\delta_{x=\frac{3}{4}}\right] +g(0) \left[\delta'_{x=-\frac{3}{4}}-\delta'_{x=\frac{3}{4}}\right].$$
It is easy to check that
$$\norm{Lg}_X \leq \norm{\Delta g}_X+\norm {g}_{\mathcal C^1_\phi}\leq C\norm{g}_Z.$$
\qed

% ================================================================================================================
% ================================================================================================================

\noindent {\bf Acknowledgements}\\
We would like to thank Prof. Joan Sol\'a-Morales at UPC (Barcelona) for his advise and guidance during the preparation of this paper. We would like also to thank Prof. Charles Fefferman (Princeton) for his useful suggestions. Moreover we would  like to thank Prof. Richard Tsai (UT Austin) for his numerical simulations, which provided nice intuitions for the analysis.\\
Support to the second author from the Centre de Recerca Matem\`atica  at the Universitat Aut\`onoma de Barcelona is also very gratefully acknowledged.
\\Both authors thank IPAM (UCLA) for the opportunity to collaborate.

%-------------------------------------------------------------------------------------------------------
%-------------------------------------------------------------------------------------------------------

\end{document}